\documentclass[article,
 amsmath,amssymb,
 reprint,%
]{revtex4-1}

\usepackage{graphicx}
\usepackage{dcolumn}
\usepackage{bm}

\usepackage[utf8]{inputenc}
\usepackage[T1]{fontenc}
\usepackage{mathptmx}
\usepackage{etoolbox}

\makeatletter
\def\@email#1#2{%
 \endgroup
 \patchcmd{\titleblock@produce}
  {\frontmatter@RRAPformat}
  {\frontmatter@RRAPformat{\produce@RRAP{*#1\href{mailto:#2}{#2}}}\frontmatter@RRAPformat}
  {}{}
}%
\makeatother

\usepackage[usenames]{color}
\usepackage{colortbl}

\begin{document}

\preprint{AIP/123-QED}

\title[Sample title]{Bifurcations of mode-locked periodic orbits in three-dimensional maps}

\author{Sishu Shankar Muni, Soumitro Banerjee}%
\affiliation{ 
Department of Physical Sciences\\
Indian Institute of Science Education and Research Kolkata\\
Campus Road, Mohanpur, West Bengal, 741246\\
India\\
{\rm{E-mail: sishushankarmuni@gmail.com}}
}%


\date{\today}

\begin{abstract}

  In this paper, we report the bifurcations of mode-locked periodic
  orbits occurring in maps of three or higher dimensions. The `torus'
  is represented by a closed loop in discrete time, which contains
  stable and unstable cycles of the same periodicity, and the unstable
  manifolds of the saddle. We investigate two types of `doubling' of such
  loops: in (a) two disjoint loops are created and the iterates toggle
  between them, and in (b) the length of the closed invariant curve is
  doubled. Our work supports the conjecture of Gardini and Sushko,
  which says that the type of bifurcation depends on the sign of the
  third eigenvalue. We also report the situation arising out of
  Neimark-Sacker bifurcation of the stable and saddle cycles, which
  creates cyclic closed invariant curves. We show interesting types of
  saddle-node connection structures, which emerge for parameter values where
  the stable fixed point has bifurcated but the saddle has not, and vice versa.
 
\end{abstract}

\maketitle
\section{Introduction}

Mode-locked periodic orbits occur when a dynamical system has two
commensurate frequencies that are not harmonically related to each
other. Such cycles always occur as stable-unstable pairs, and are
located on a closed invariant curve formed by the unstable manifolds
of the saddle. Such closed invariant curves occurring in maps
represent tori in continuous time. In this paper we investigate
bifurcations of such `resonant tori' in three-dimensional maps.

Much of the earlier work on the bifurcation of tori focused on ergodic
tori or quasiperiodic orbits.  The doubling of quasiperiodic orbits
was first reported by Kaneko \cite{Kaneko83} and Arneodo et
al. \cite{Spiegel83} in the dynamics of some three- and
four-dimensional maps. Since then, such bifurcations have been
reported in maps resulting from the discretization of continuous-time
systems such as impulsive Goodwin's oscillator \cite{Zhu21},
3-dimensional Lotka-Volterra model \cite{Banerjee12,Gardini87},
vibro-impacting systems \cite{ding2004interaction}, four coupled
oscillators \cite{Ash95}, radiophysical systems \cite{Anish05},
etc. Sekikawa et al. demonstrated a sequence of length-doublings
\cite{Seki01,Seki21} finally resulting in chaos.

  Torus doubling critical point is a special point in the parameter
  space where the regions of occurrence of torus, doubled torus, and
  strange non-chaotic attractor meet. Such critical points were
  found in a nonlinear electronic circuit with quasiperiodic drive
  \cite{Seleznev2}. Scaling laws applicable in the vicinity of a
  critical point were obtained for a forced logistic map
  \cite{forcedlogistic}.

  For a better understanding on the bifurcations of invariant tori,
  iterative schemes for computation and continuation of
  two-dimensional stable and unstable manifolds of mode-locked
  orbits were proposed in \cite{Meiss06,Vitolo08}. In piecewise-smooth
  maps, the birth of a bilayered torus via border-collision bifurcation was
  reported in \cite{Zhu06} and bifurcation of such a torus via
  heteroclinic tangencies was discussed.  Recently, a review on the
  dynamical consequences of torus doublings resulting in the formation
  of Shilnikov attractors were discussed in \cite{Gon21}.

  After Grebogi, Ott and Yorke \cite{grebogi1983three} showed that
  three-frequency quasiperiodicity can be stable over a parameter range,
  its occurrence, observed in the form of a torus in the discrete-time
  phase space, has been reported in many systems including symmetrical
  rigid tops \cite{nstoro6}, vibro-impacting systems \cite{Vibr15},
  four coupled Chua's circuits \cite{zhong1998torus}, a simple autonomous
  system \cite{kuznetsov2016simplest}, and a ladder of four coupled
  van der Pol oscillators \cite{Kom16}.

  The merger and disappearance of a stable and an unstable closed
  invariant curves was found to be responsible for the creation of
  strange nonchaotic attractors
  \cite{Hammel94,FEU95,SNAbook06}. Behaviors close to tori-collision
  terminal points were studied in \cite{kuznetsov2000critical}.

\begin{figure}[b]
  \centering
  \includegraphics[width=2.5in]{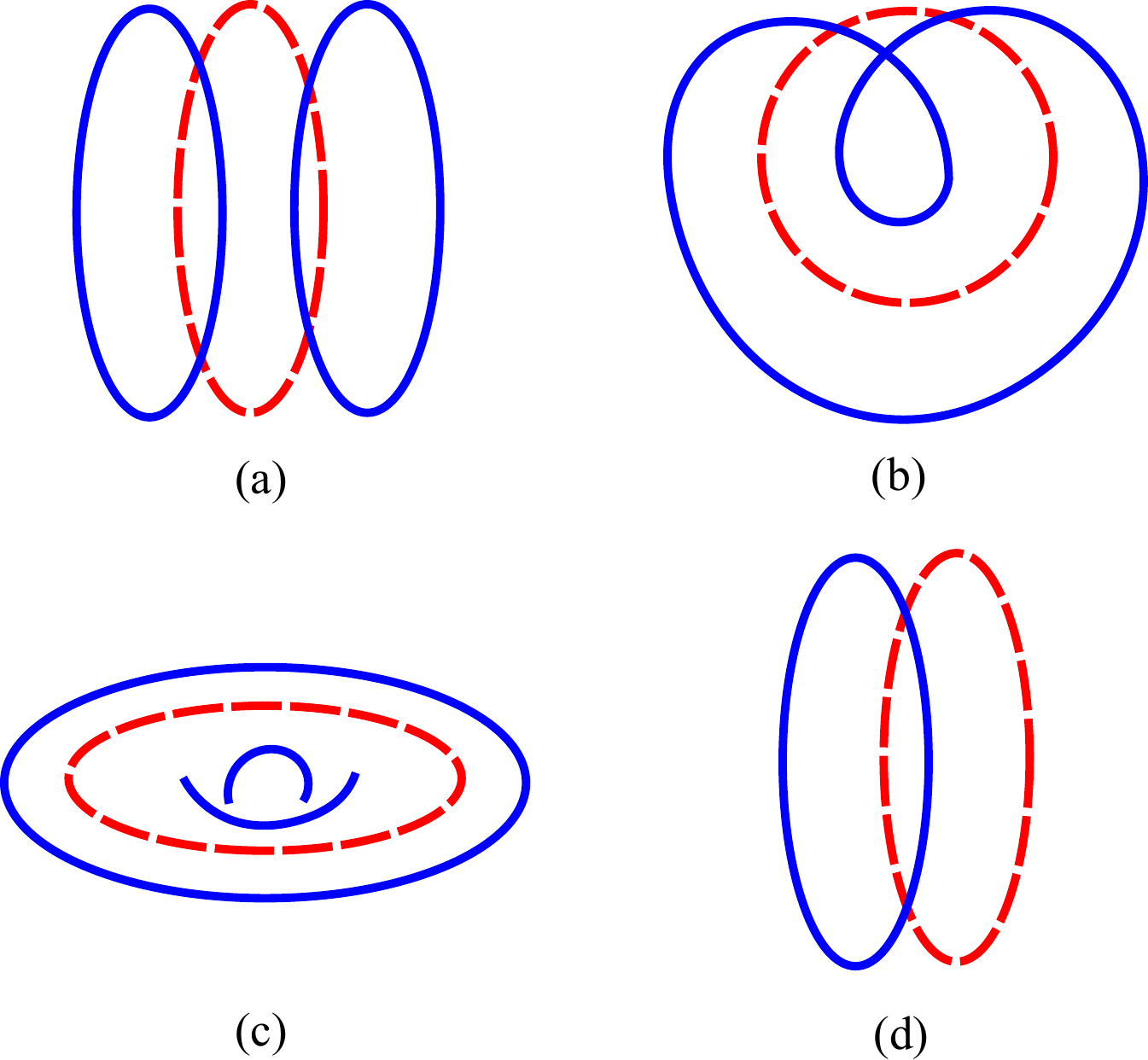}

 \caption{The four types of bifurcations of a quasiperiodic orbit \cite{Banerjee12}. Unstable objects are shown with red dashed lines and stable objects are marked in blue.\label{quasifig}}
\end{figure}

These observations have been summarized in \cite{Banerjee12}, by
showing that a stable quasiperiodic orbit can bifurcate in four
possible ways. These are illustrated in
Fig.~\ref{quasifig}. Fig.~\ref{quasifig}(a) and (b) show two ways of
`doubling': in (a) two disjoint loops are created and the iterates
toggle between them, and in (b) the length of the closed invariant
curve is doubled. The bifurcation shown in Fig.~\ref{quasifig}(c)
leads to a torus in discrete time which represents three-frequency
quasiperiodicity. Fig.~\ref{quasifig}(d) shows a situation where a
stable and an unstable tori merge and disappear (or are created if the
parameter is varied in the opposite direction). Banerjee et al. explained these bifurcations based on the method of `second Poincar\'e section' \cite{Banerjee12}.

The above lines of work concern bifurcations of quasiperiodic orbits.
However, mode-locked periodic orbits represent a generic
case. Such orbits possess rational rotational numbers and
appear as shrimp shaped regions in two-dimensional parameter space
\cite{Ga94} that cover a larger range of parameters than the
quasiperiodic orbit.  Therefore, in this paper we explore the
ways in which ‌mode-locked periodic orbits may bifurcate.

Although much research has been done on the doubling of quasiperiodic
orbits, much less research attention has been devoted to the doubling
of resonant tori or mode-locked periodic orbits. Gardini and Sushko
\cite{Gardini12} proposed a conjecture regarding the conditions
for different types of doubling bifurcations of closed invariant curves
related to mode-locked periodic orbits. However, to date these have
not been tested using systems that exhibit these bifurcations. In this
paper we fill that gap and validate their conjecture.

We also explore the birth of a third frequency in a mode-locked
periodic orbit, which results in the occurrence of more than two
closed invariant curves in the phase space and the iterates cyclically
move among them.

Both types of torus doubling are related to the occurrence of period
doubling bifurcation and the birth of multiple loops is related to the
occurrence of a Neimark-Sacker bifurcation in a fixed point. Since a
saddle cycle as well as a stable cycle occur on the closed invariant
curve, an interesting situation emerges: it is possible that one of
them undergoes a bifurcation but the other does not. Such situations
lead to atypical structures of the closed invariant curve, which are
also reported in this paper.

The paper is organized as follows: In \S
\ref{sec:QualTheoryModeLocked}, we recapitulate the conjecture and
discuss the mechanism behind the doubling of mode-locked orbits in
three-dimensional maps. In \S \ref{sec:DisjointLoop}, we give an
example of the doubling phenomenon of a mode-locked periodic orbit in
a three-dimensional smooth map where two disjoint loops are formed. In
\S \ref{sec:LengthDoubled}, we provide example of the case where
the length of manifold doubles and they lie on a non-orientable
M\"{o}bius strip.  In \S \ref{sec:collision}, collision of two closed
invariant curves formed out of a saddle-node connection and a
saddle-saddle connection is discussed. \S \ref{sec:Conclusions}
presents the conclusions and future research questions arising out of
this paper.

\section{Qualitative theory of doubling of mode-locked orbit}
\label{sec:QualTheoryModeLocked}

Closed invariant curves are born through Neimark-Sacker bifurcation of a stable fixed point. This can lead to a few scenarios:

\begin{enumerate}
\item A supercritical Neimark-Sacker bifurcation leading to the formation of a stable closed invariant curve.
  \begin{enumerate}
  \item There are a dense set of points on the curve, which implies an irrational frequency ratio. This leads to the formation of a quasiperiodic orbit.\label{quasi}
    \item There are a finite number of stable and saddle points and the invariant closed curve is composed of a union of these points and the unstable manifolds of the saddle points. This leads to the formation of a mode-locked periodic orbit. 
  \end{enumerate}

\item A subcritical Neimark-Sacker bifurcation leading to an unstable closed invariant curve.
  \begin{enumerate}
  \item The unstable closed invariant curve may be a quasiperiodic orbit;
    \item The unstable closed invariant curve may be a mode-locked periodic orbit.
  \end{enumerate}

\end{enumerate}

In this paper we consider the bifurcations of a closed invariant curve related to Case 1(b).

\begin{figure}[tbh]
    \centering
    \includegraphics[width=0.8\columnwidth]{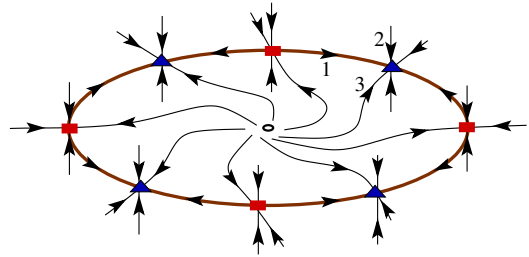}
\caption{Schematic structure of a closed invariant curve for a mode-locked periodic orbit. The stable cycle is marked by blue triangles and the saddles by red rectangles. The directions associated with eigenvalues 1, 2 and 3 of the stable cycle are shown for one of the points.}\label{invariant}
  \end{figure}

Let us consider a three-dimensional map with a stable cycle (say
$C^{N}$) and saddle cycle (say $C^{S}$). The union of these points and
the unstable manifolds of the saddle cycles form the closed invariant
curve (see Fig.~\ref{invariant}). Let us denote the eigenvalues of the stable periodic orbit as
$\lambda_{i}^{N}$, $i=1,2,3$ and the eigenvalues of the saddle
periodic orbit as $\lambda_{i}^{S}$, $i=1,2,3$.

In order for the closed invariant curve to exist, one of the three
eigenvalues must be positive. Let the positive eigenvalues be denoted
by $\lambda_{1}^{N}, \lambda_{1}^{S}$ for the node and saddle periodic
points, respectively. The corresponding stable manifolds of the nodes
and the unstable manifolds of the saddles form the saddle-node
connections on the closed invariant curve.

Following \cite{Gardini12}, let us first consider the case of doubling of the closed invariant
curve. In order for a doubling to occur, there must be a flip
bifurcation, i.e., one of the eigenvalues should pass through
$-1$. Let this eigenvalue be denoted as $\lambda_{2}^{N}$ for the node
and as $\lambda_{2}^{S}$ for the saddle. These eigenvalues, associated
with the saddle and the node, need not pass through $-1$ at the same
parameter value. Indeed, the saddle and the node can undergo flip
bifurcation at different parameter values. We will see this in two
explicit examples considered later in this paper.

\begin{figure*}[htbp]
	\centering
\includegraphics[width=7cm]{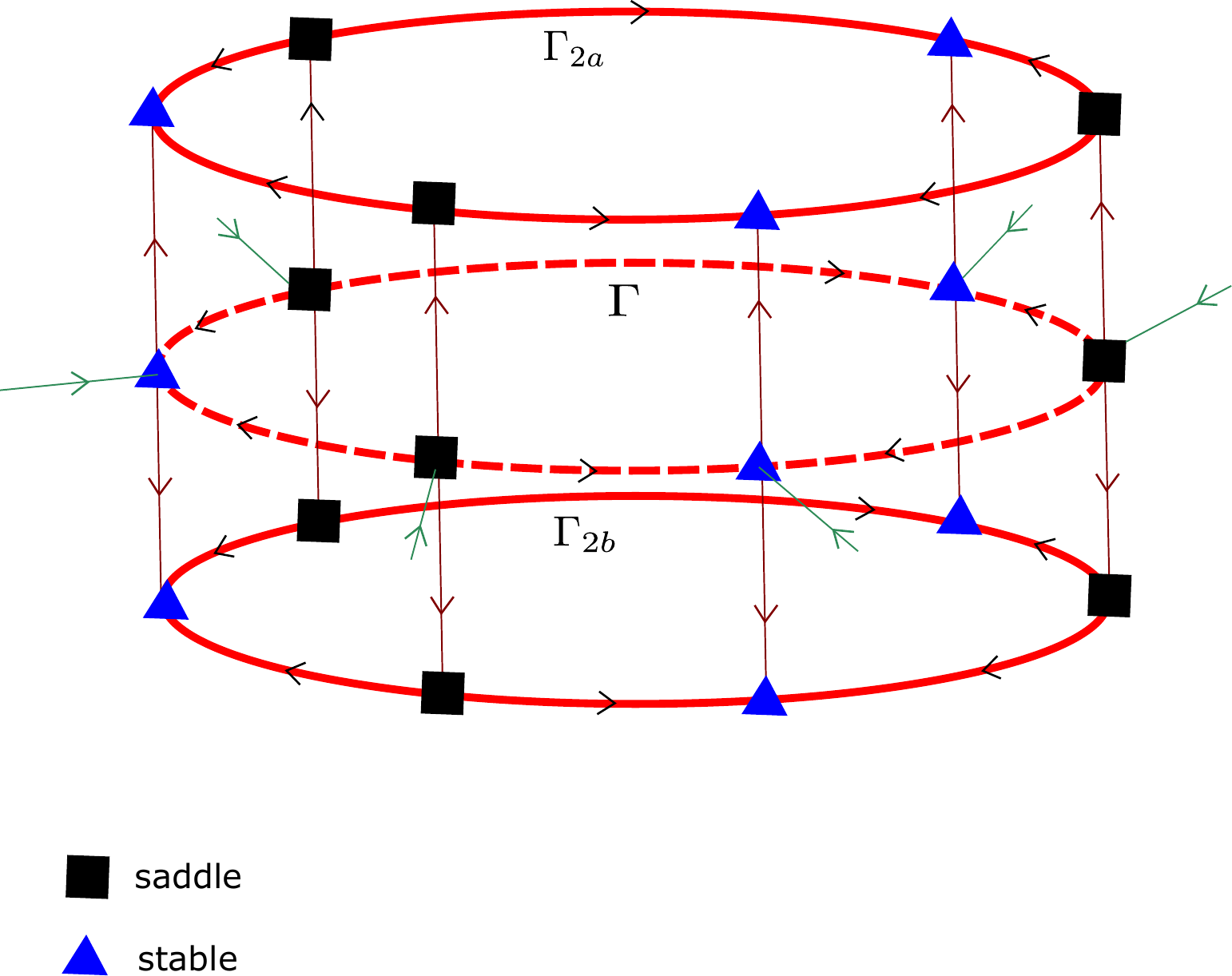}\hspace{1cm}
 \includegraphics[width=6cm]{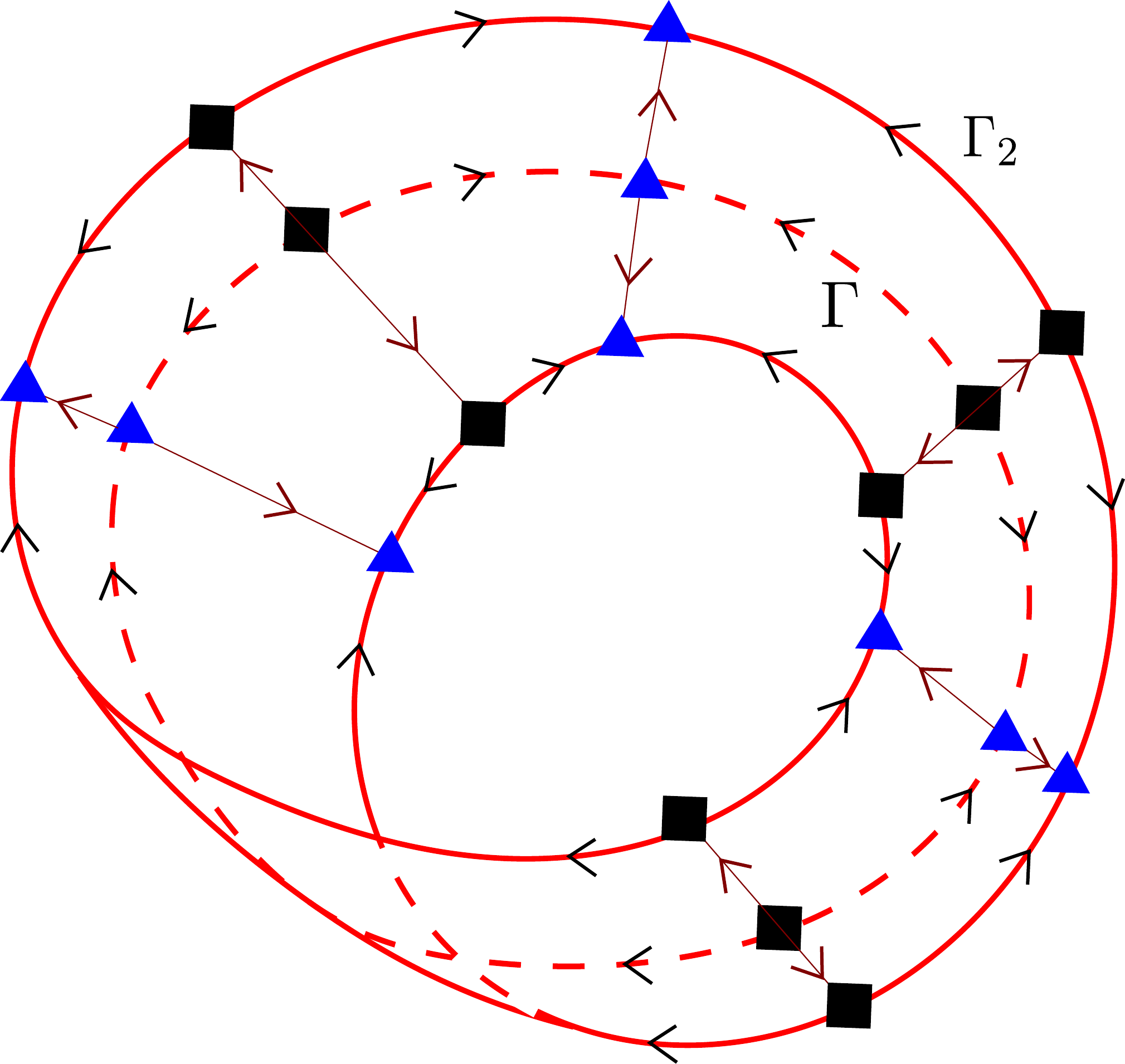} 

{\small \bf (a)} \hspace{2.6in}
{\small \bf (b)}
\caption{In (a), the doubling bifurcation related to the positive third eigenvalue leading to disjoint loops. In (b),  the doubling bifurcation related to negative third eigenvalue and the mode locked periodic orbit lies on the surface of a M\"{o}bius strip.}
\label{fig:ManifoldSketch}
\end{figure*}

Consider a situation where both the saddle $C^{S}$ and the node
$C^{N}$ have undergone flip bifurcation (see
Fig.~\ref{fig:ManifoldSketch}).  A doubled stable orbit (say $C^{2N}$)
has emanated from $C^{N}$, which has turned into a flip-saddle with
 one unstable direction.  A similar scenario occurs with the
saddle periodic point $C^{S}$, which doubles by making another of its
branches unstable (two independent unstable directions). Two saddle
points (say $C^{2S}$) emanate on both sides of the previous saddle
($C^{S}$).
Two distinctly different stable behaviors can result from this bifurcation.
Gardini and Sushko \cite{Gardini12} conjectured that the shape of the resulting manifold is determined by the third eigenvalue $\lambda_{3}^{N}, \lambda_{3}^{S}$. The sign of the third eigenvalue classifies the bifurcation into two types.

\subsection{Positive third eigenvalue ($\lambda_{3}^{N} > 0, \lambda_{3}^{S} > 0$)}

If the third eigenvalue associated is positive, then the trajectories
will be converging on one side of the manifold and will give a
geometric structure like Fig.~\ref{fig:ManifoldSketch}(a). This leads
to the formation of two disjoint cycles $\Gamma_{2a}, \Gamma_{2b}$. We
note that the two disjoint cycles are cyclically invariant, i.e.,
the trajectory on the periodic orbits toggle between the upper cycle
$\Gamma_{2a}$ and the lower cycle $\Gamma_{2b}$. The two disjoint closed loops
are invariant in the sense that they act as a single loop for the
second iterate of the map.

The two attracting disjoint cycles bound a strip or a manifold whose
shape is topologically a cylinder. The manifold (topological cylinder)
consists of the unstable closed cycle $\Gamma$ and the two stable
disjoint cycles $\Gamma_{2a}$, and $\Gamma_{2b}$. The cylinder
manifold is orientable.   Each of the node and saddle periodic
points ($C^{N}, C^{S}, C^{2N}, C^{2S}$) have three independent
eigen-directions, which are
represented in Fig.~\ref{fig:ManifoldSketch} with different colours.

\subsection{Negative third eigenvalue ($\lambda_{3}^{N} < 0, \lambda_{3}^{S} < 0$)}

If the third eigenvalue is negative, the shape of the invariant set is
a M\"{o}bius strip (a non-orientable manifold). The boundary of the
M\"{o}bius strip constitutes the
saddle-node connection of the doubled saddle and node mode-locked
periodic orbits. This is qualitatively shown in
Fig.~\ref{fig:ManifoldSketch}(a).

We now illustrate the above bifurcations with a few examples.

\section{Example: Disjoint two-loop mode-locked orbit}
\label{sec:DisjointLoop}

In this section we consider the Mira map studied in \cite{Kara21}, which is a three-dimensional map given by 
\begin{equation}
    \begin{aligned}
    x' &= y,\\
    y' &= z,\\
    z' &= Bx + Cy + Az - y^2,\\
    \end{aligned}
    \label{eq: ShilnikovMap}
\end{equation}

The Jacobian of the map \eqref{eq: ShilnikovMap} is
\begin{equation}
J = \begin{bmatrix}
0 & 1 & 0\\
0 & 0 & 1\\
B & C -2y & A 
\end{bmatrix},
\label{eq:JacobianShilnikov}
\end{equation}
whose determinant is det$(J) = B$. For $B<0$, the map is
orientation-reversing everywhere, and for $B>0$, the map is
orientation-preserving everywhere. We have observed disjoint two-loop  mode-locked periodic orbit in the
orientation-reversing regime. In \cite{Kom16}, authors have mentioned
that such disjoint doublings of quasiperiodic orbit can be found in
orientation-preserving systems as well but in higher dimensional
systems of dimensions greater than or equal to four.
 
From the one-parameter bifurcation diagram with the variation of
parameter $B$, see Fig. \ref{fig:OneParamBif_Shilnikov_B}(a), we
observe that a fixed point undergoes a Hopf bifurcation to a
quasiperiodic orbit and then through a saddle-node bifurcation, a
mode-locked period-five orbit is formed. It then doubles to a
period-10 mode-locked orbit. With further increase in the parameter
$B$, we observe subsequent doublings of mode-locked orbit, after which
it transits to chaos.

\begin{figure}[tbh]
\centering
\setlength{\unitlength}{1cm}
\includegraphics[width=0.48\columnwidth]{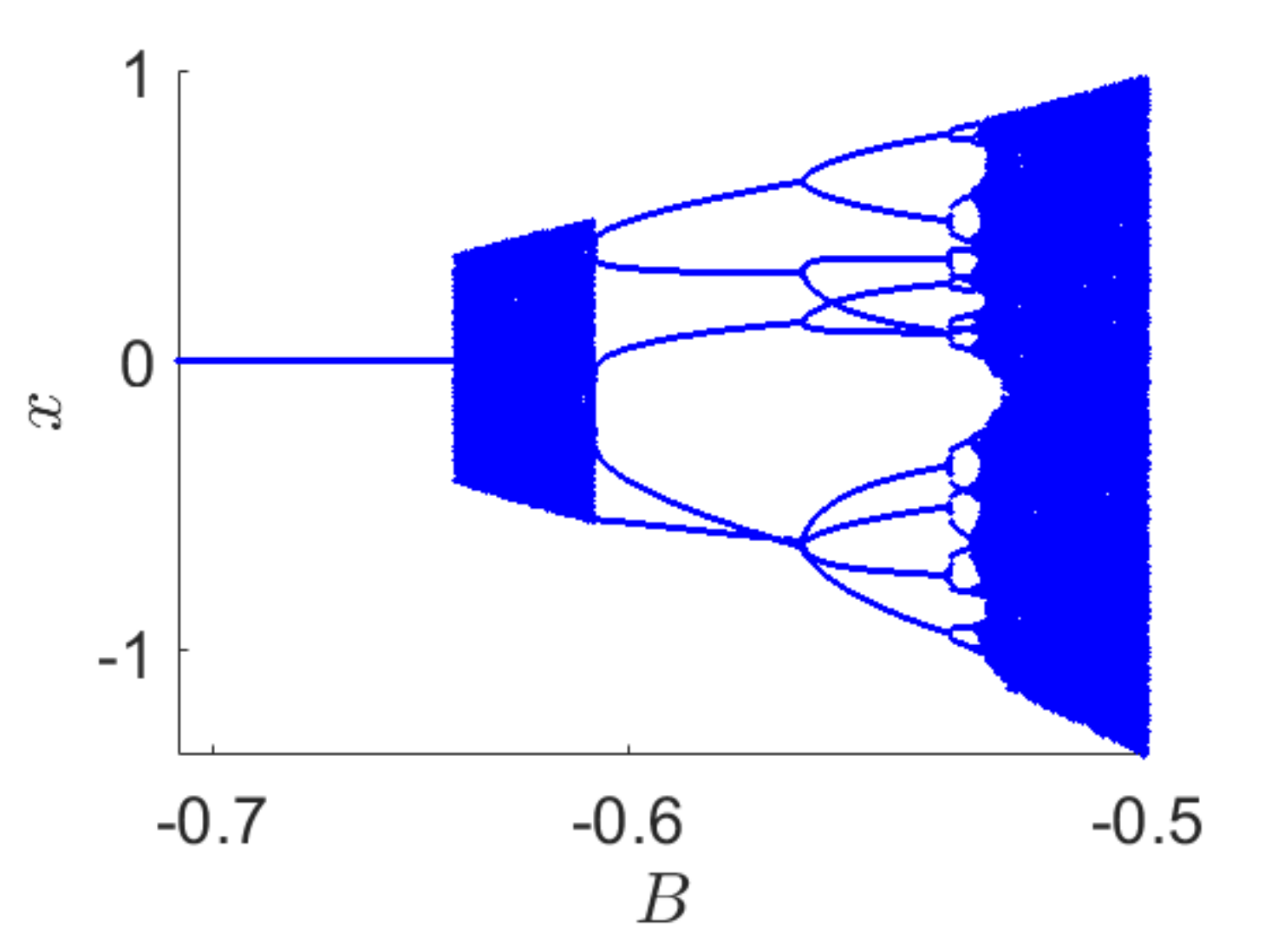}
\includegraphics[width=0.48\columnwidth]{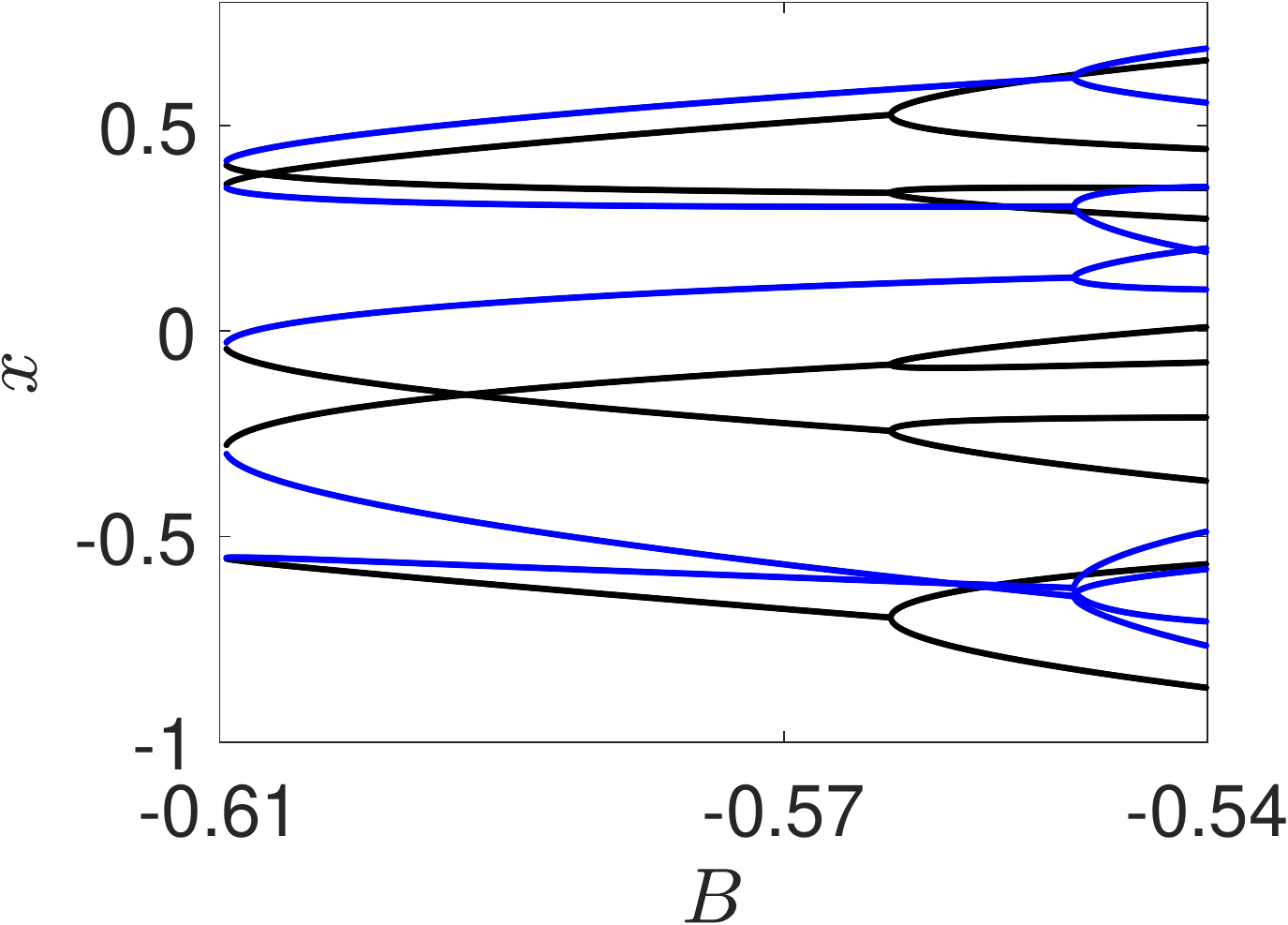}

{\small  (a)} \hspace{1.5in}
{\small  (b)}

\caption{(a) One-parameter bifurcation diagram of $x$ vs $B$. The parameters are $A=-2.269, C =-2.1$. At $B =-0.58$, we observe a period-5 mode-locked orbit. At $B=-0.54$, we observe a period-10 mode-locked orbit. When $B$ is further increased, it undergoes subsequent period-doubling bifurcations leading to chaos. (b)  One-parameter bifurcation diagram showing
  both the saddle (black) and stable (blue) periodic points. 
 }
\label{fig:OneParamBif_Shilnikov_B}
\end{figure}

The doubling of the mode-locked orbit can be better understood when we
continue the saddle and stable mode-locked periodic orbits with respect
to the parameter
$B$. We have used the multidimensional Newton-Raphson method to locate the
saddle cycle. Fig. \ref{fig:OneParamBif_Shilnikov_B}(b) shows that
the saddle doubles first, followed by the doubling of the node.  The
computation of the eigenvalues over the same parameter range of $a$
(see Fig. \ref{fig:Eigenvaue_Saddle_Shilnikov_B}) reveals that, as the
parameter $B$ increases, the second eigenvalue reaches $-1$ at $B=-0.5627$. In
Fig. \ref{fig:Eigenvaue_Saddle_Shilnikov_B}(b), we see that the
second eigenvalue of the node reaches $-1$ at $B=-0.55$ and hence the
stable period-five orbit bifurcates to a period-ten orbit.


\begin{figure}[!htbp]
\centering
\includegraphics[width=0.48\columnwidth]{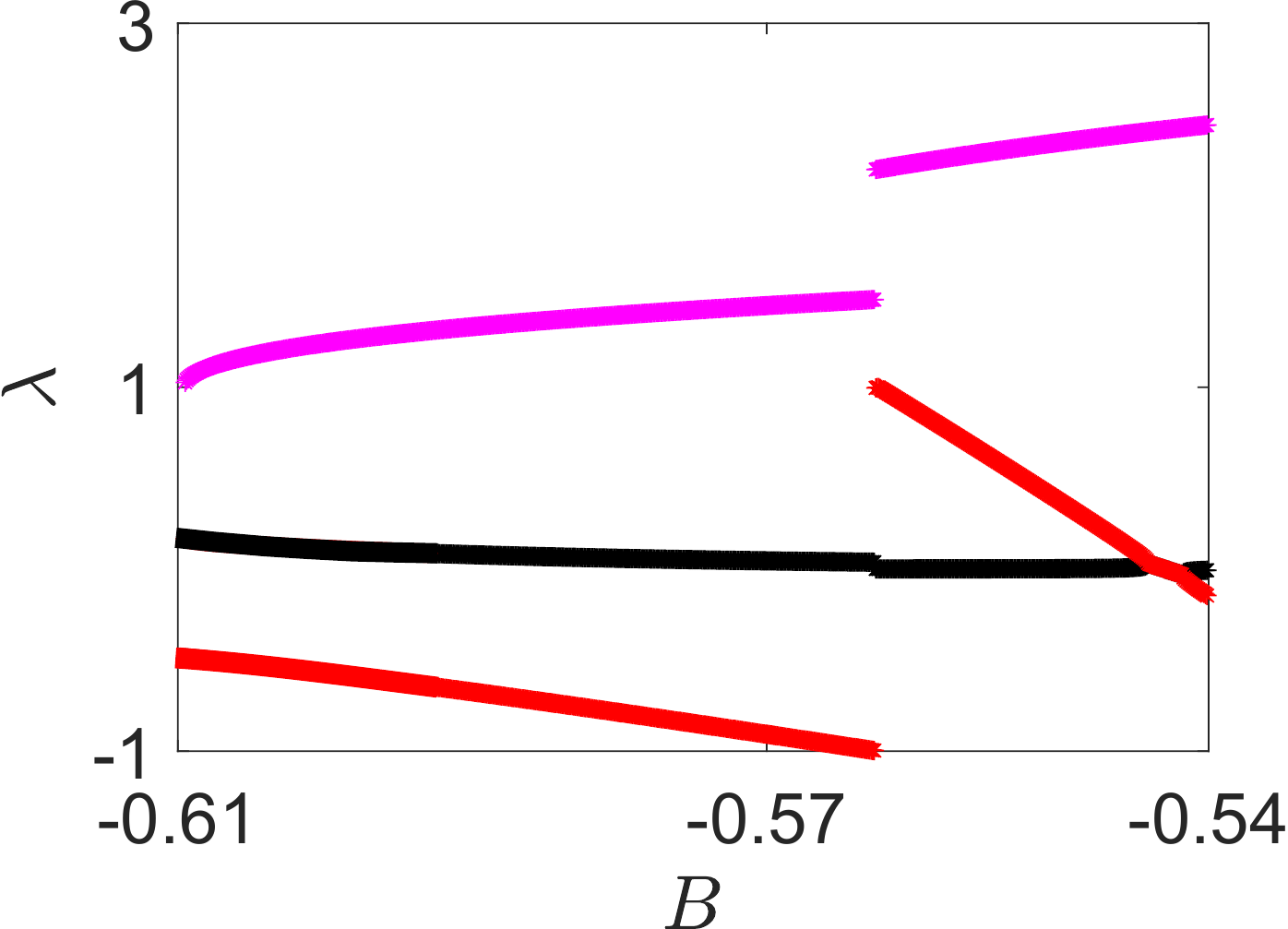}\hfill
\includegraphics[width=0.48\columnwidth]{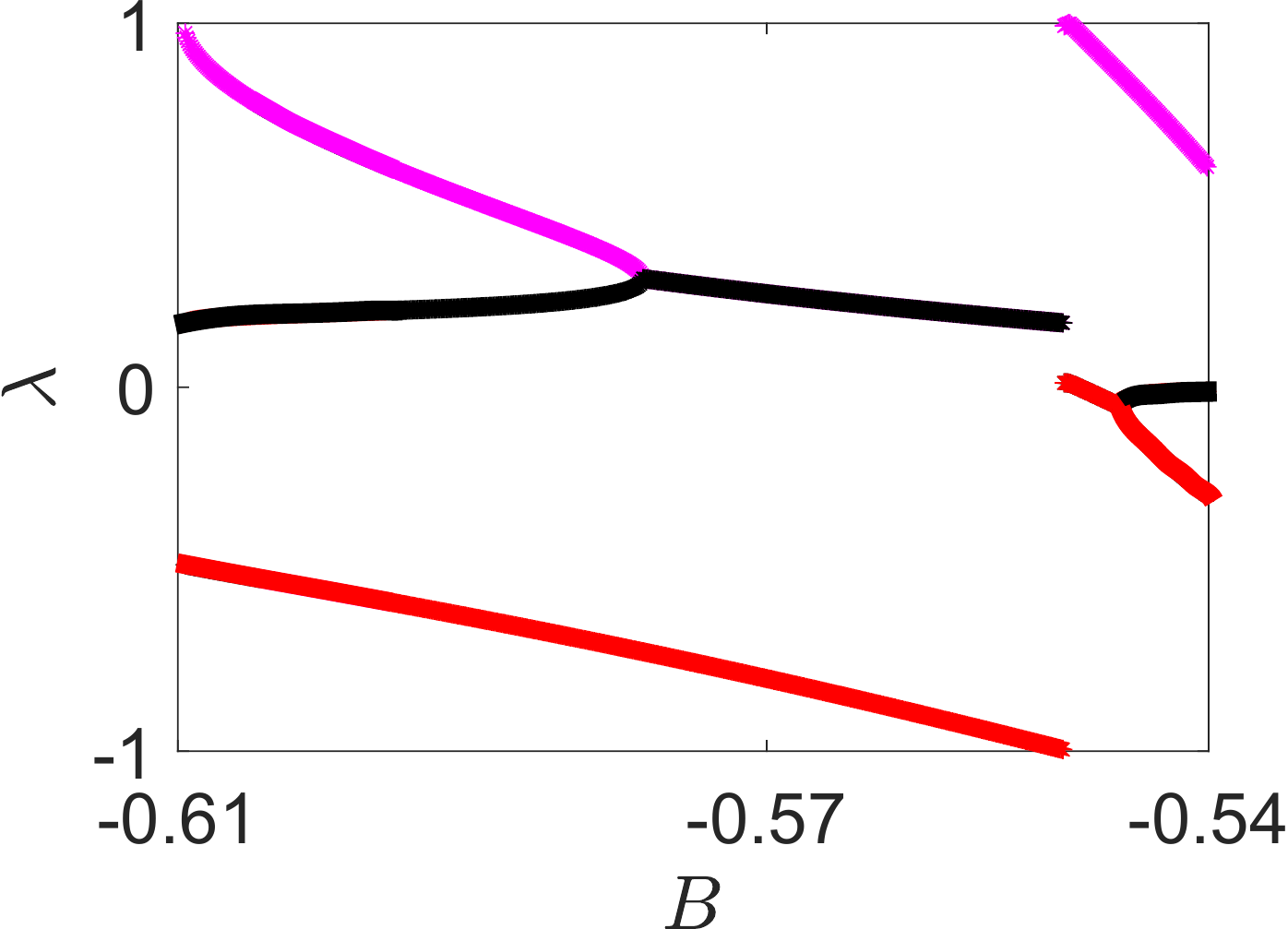}\\
{\small  (a)} \hspace{1.6in}
{\small  (b)}
\caption{ Continuation of the eigenvalues as the parameter $B$
  is varied: (a) for the saddle cycle and (b) for the stable cycle. Three eigenvalues are differentiated by different colours: eigenvalue $\lambda_1$ is denoted by magenta,
    $\lambda_2$ by red, and $\lambda_3$ by black. }
\label{fig:Eigenvaue_Saddle_Shilnikov_B}
\end{figure}

\begin{figure}[htb]
\centering
\includegraphics[width=0.49\columnwidth]{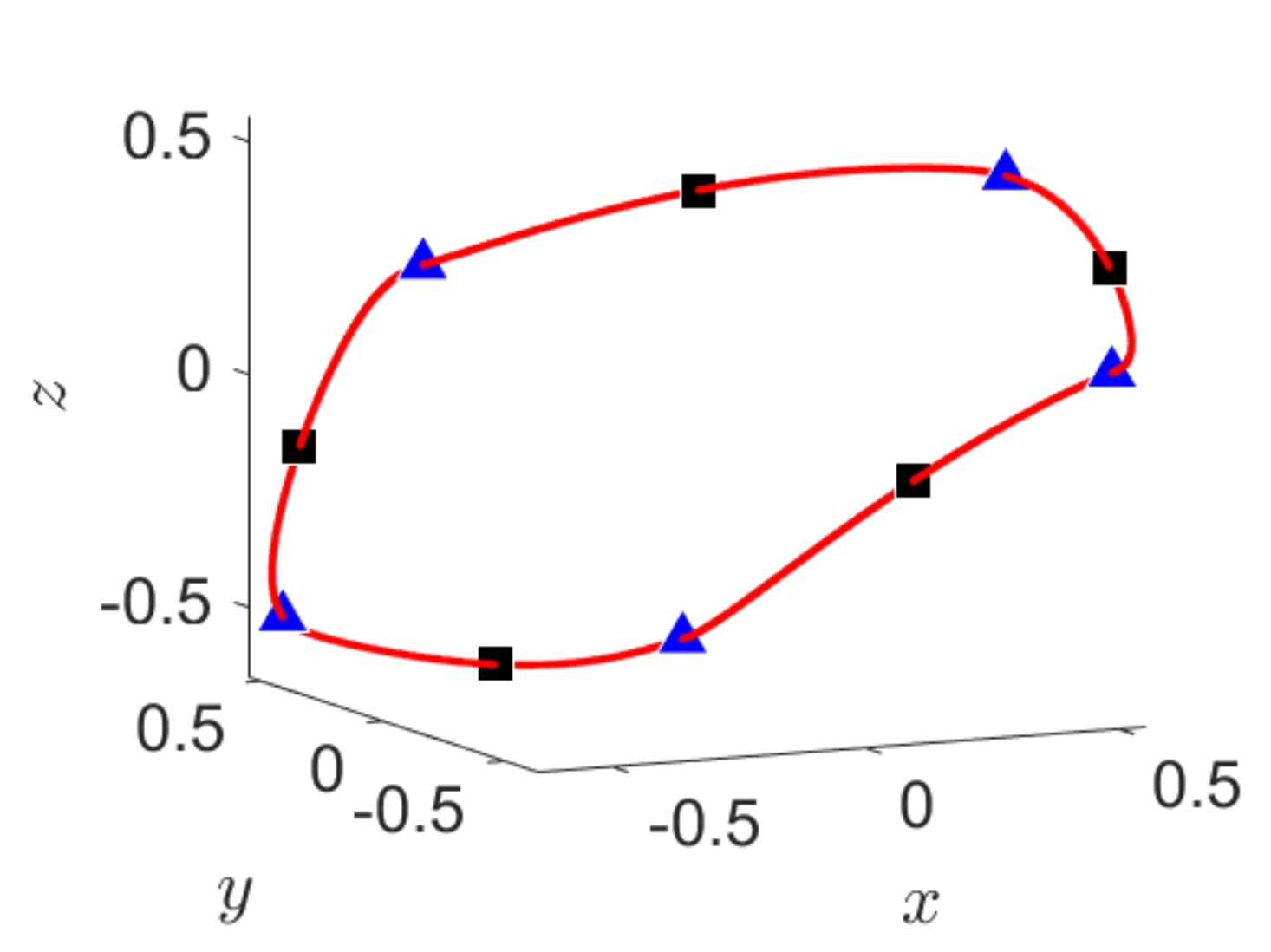}\hfill
\includegraphics[width=0.49\columnwidth]{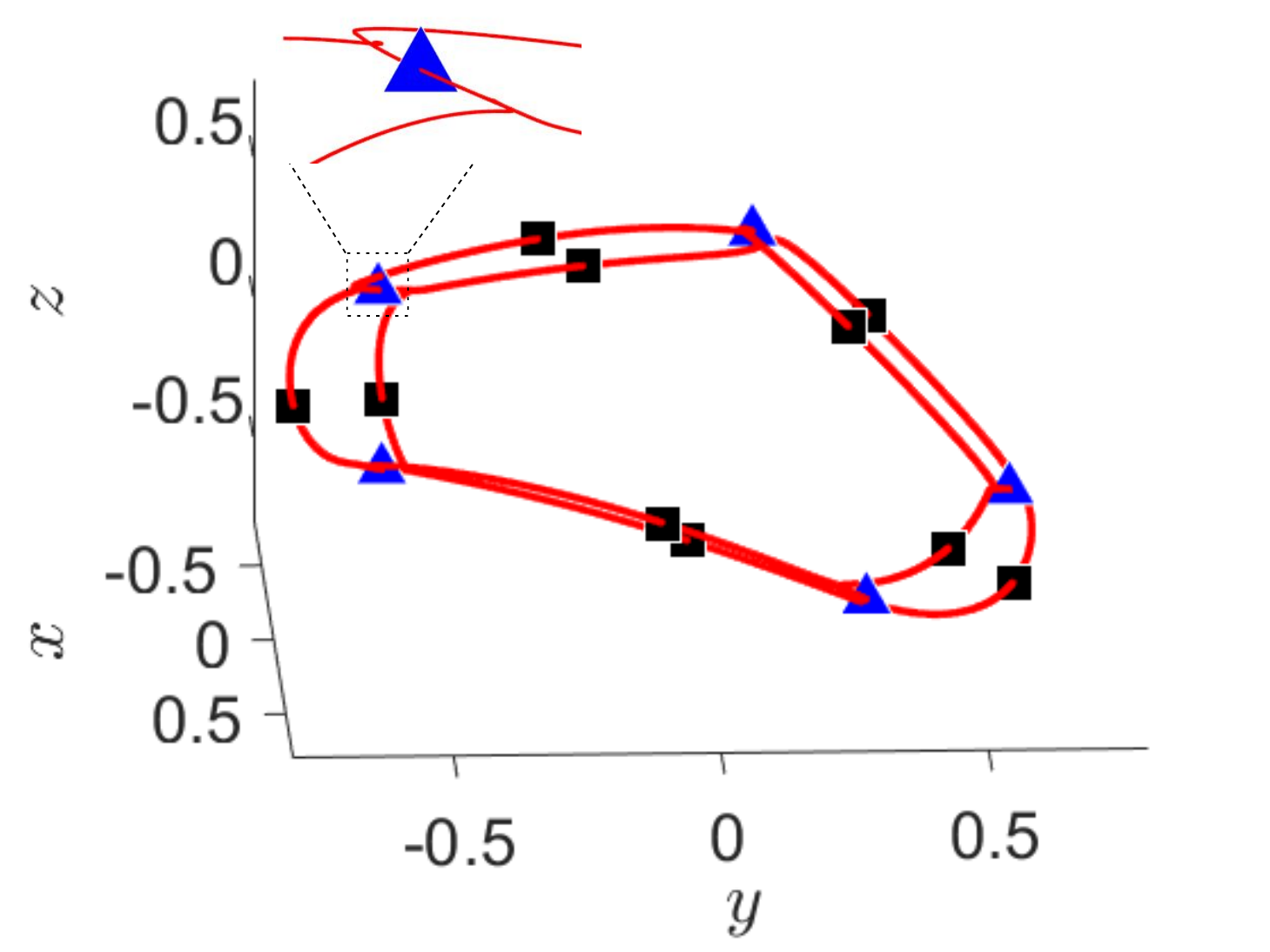}\\
{\footnotesize  (a)} \hspace{1.6in}
{\footnotesize  (b)}\\
\includegraphics[width=0.49\columnwidth]{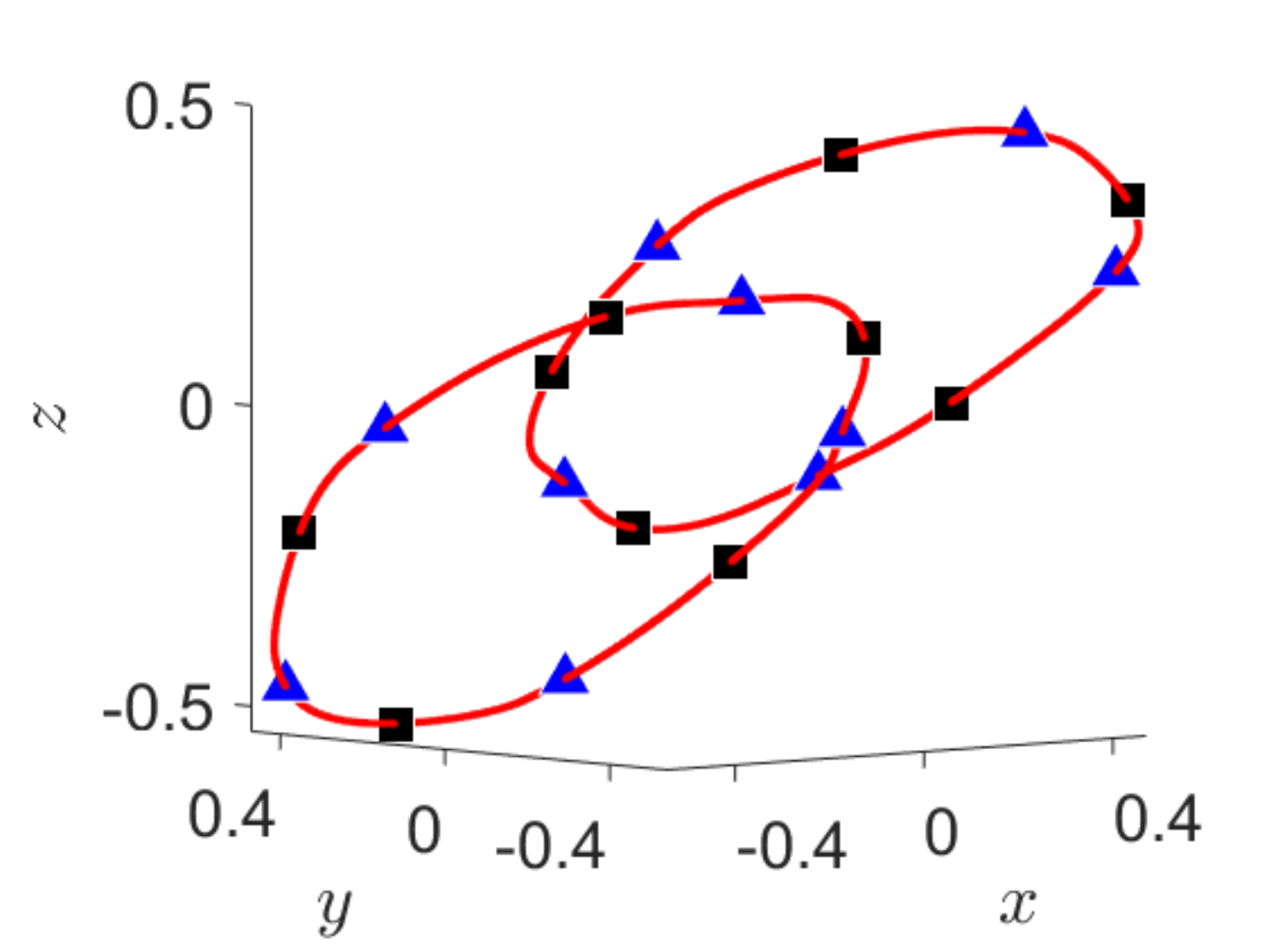}
{\footnotesize (c)}
\caption{(a) The closed invariant curve formed by the
  unstable manifolds (in red) of the saddle period-5 periodic orbit at
  $B=-0.58$. The stable mode-locked period-5 orbit is denoted by blue
  triangles and the saddle period-5 orbit is denoted by black
  squares. (b) The structure of the closed invariant curve for
  $B = -0.555$, when the saddle has doubled but the stable cycle has not yet doubled. The zoomed portion shows the connection through the stable fixed point. (c) At $B =-0.54$, two disjoint loops are
  formed by the unstable manifolds of two disjoint yet cyclic saddle
  period-10 mode-locked orbits. The other parameters are $ A = -2.269$ and $C = -2.1$.}
\label{fig:SingleLoopShilnikov}
\end{figure}

At $B = -0.58$, a stable period-5 mode-locked orbit exists along with
a period-5 saddle.  We compute the unstable manifolds of the saddle
cycle using the method of fundamental domains
\cite{MuMcSi21}, which form a saddle-node
connection. The resulting invariant closed curve is shown in Fig. \ref{fig:SingleLoopShilnikov}(a).

As we increase the parameter $B$, a stage is reached where the saddle
has doubled but the stable periodic orbit has not. At
$B=-0.555$, the saddle period-10 orbit coexists with a stable
period-five orbit. The manifolds of the saddle point are shown in Fig.~\ref{fig:SingleLoopShilnikov}(b). We notice a complex structure in which two loops have formed, but these are joined by branches that connect with the stable period-5 orbit.

Next, we consider a parameter where both the saddle and stable
mode-locked periodic orbit have doubled. At $B=-0.54$, we see a mode
locked period-10 orbit and the unstable manifolds of the
period-doubled saddle cycle shows two disjoint loops (Fig.~\ref{fig:SingleLoopShilnikov}(b)).

Note that this disjoint mode-locked periodic orbit does not imply
bistability, rather the periodic points of the disjoint loops are
cyclically visited. The second iterate of the map shows a single loop.

The eigenvalues of the saddle and stable periodic points before the
bifurcation are shown in Table~\ref{eigens2}.  We note that before the
torus doubling bifurcation, $\lambda_1$ is positive, $\lambda_2$ is
negative (which subsequently crosses $-1$), and $\lambda_3$ is
positive.  Thus, our observation supports the conjecture by Gardini
and Sushko.

\begin{table}[htbp]
\centering
\begin{tabular}{ |c|c|c|c| } 
 \hline
 Type & $\lambda_{1}$ & $\lambda_{2}$ & $\lambda_{3}$\\ 
 \hline
 Stable  & 0.3550 &-0.7131 & 0.2593\\
 \hline
 Saddle & 1.3963 & -0.7878 & 0.0597\\ 
 \hline
\end{tabular}
  \hspace{0.5in}

  \caption{The eigenvalues of stable and saddle period-5 cycles for $B=-0.58$. 
    \label{eigens2}}
\end{table}

\section{Example: Length doubled mode-locked orbit}
\label{sec:LengthDoubled}

\begin{figure}[b]
\centering
\includegraphics[width=0.49\columnwidth]{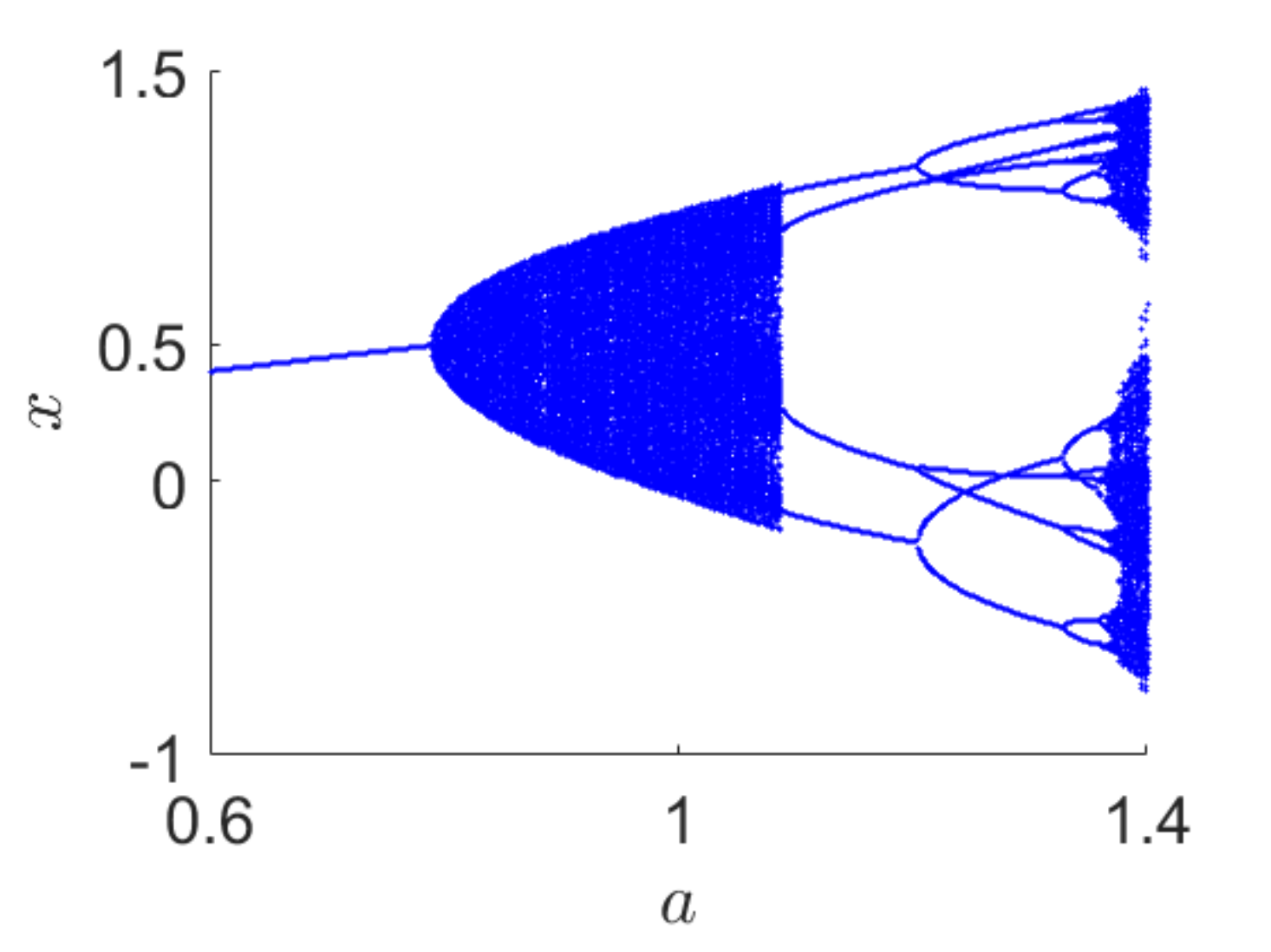}\hfill
\includegraphics[width=0.49\columnwidth]{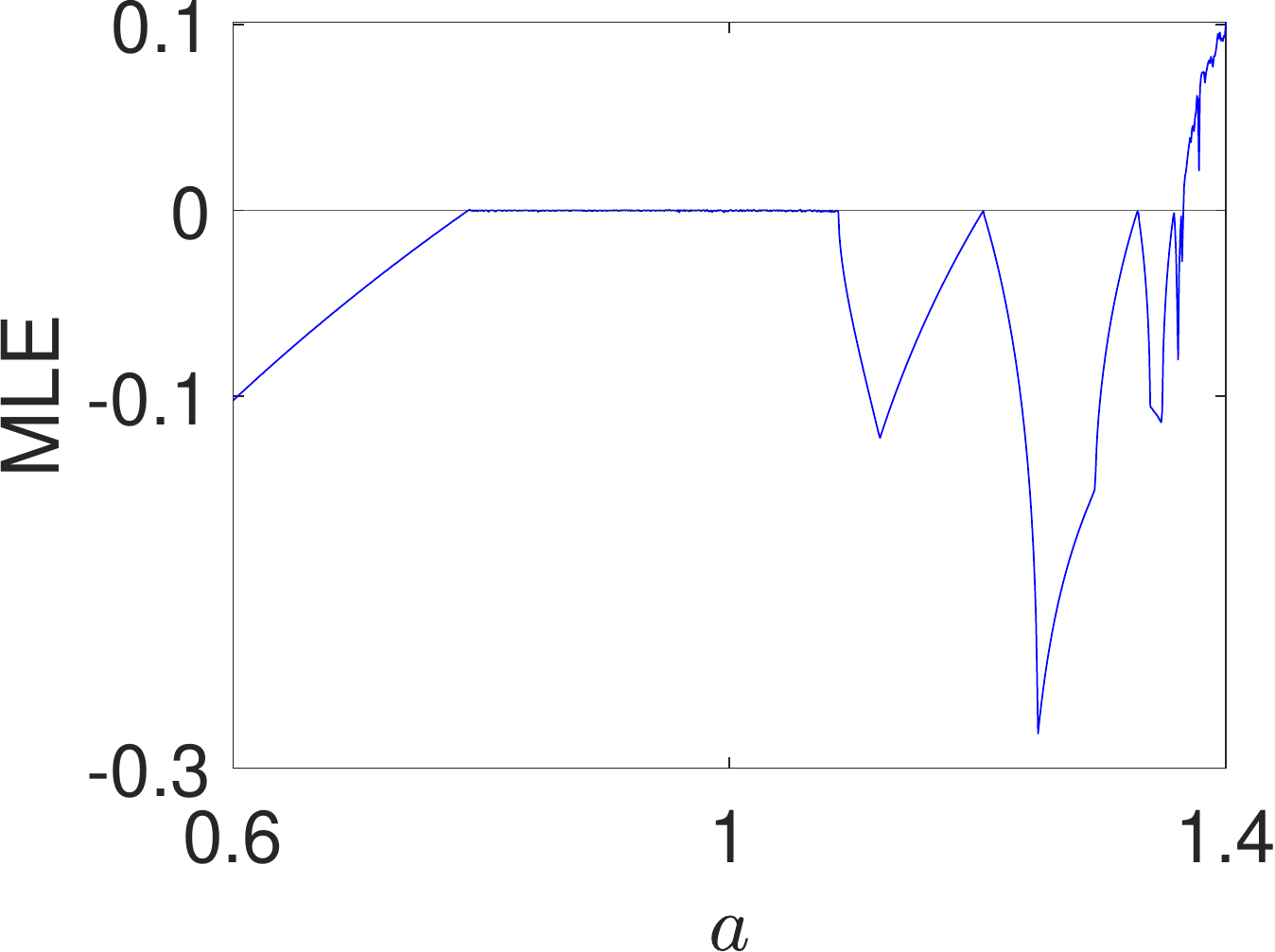}\\
{\footnotesize  (a)} \hspace{1.6in}
{\footnotesize  (b)}
\caption{In (a), a one-parameter bifurcation diagram of the generalised H\'enon map upon varying parameter $a$. In (b), plot of the maximal Lyapunov exponent in the same parameter range. The other parameter is fixed as $b=0.1$.}
\label{fig:BifurGenHenon}
\end{figure}

We consider the three-dimensional generalised H\'{e}non map \cite{Richter02} given by
\begin{equation}
    \begin{aligned}
        x' &= a - y^2 - bz,\\
        y' &= x,\\
        z' &= y,
    \end{aligned}
    \label{eq:GenHenon}
\end{equation}
where $a,b$ are the parameters. 
Fig. \ref{fig:BifurGenHenon}(a) shows a one-parameter bifurcation
diagram for this system considering $a$ as the parameter, with $b$ fixed at 0.1. A period-1 orbit bifurcates to a quasiperiodic orbit at $a\approx 0.78$ as evidenced by the maximal Lyapunov exponent being zero (Fig. \ref{fig:BifurGenHenon}(b)). A stable period-4 orbit emerges at $a \approx 1.08$.

\begin{figure}[t]
\centering
\includegraphics[width=0.6\columnwidth]{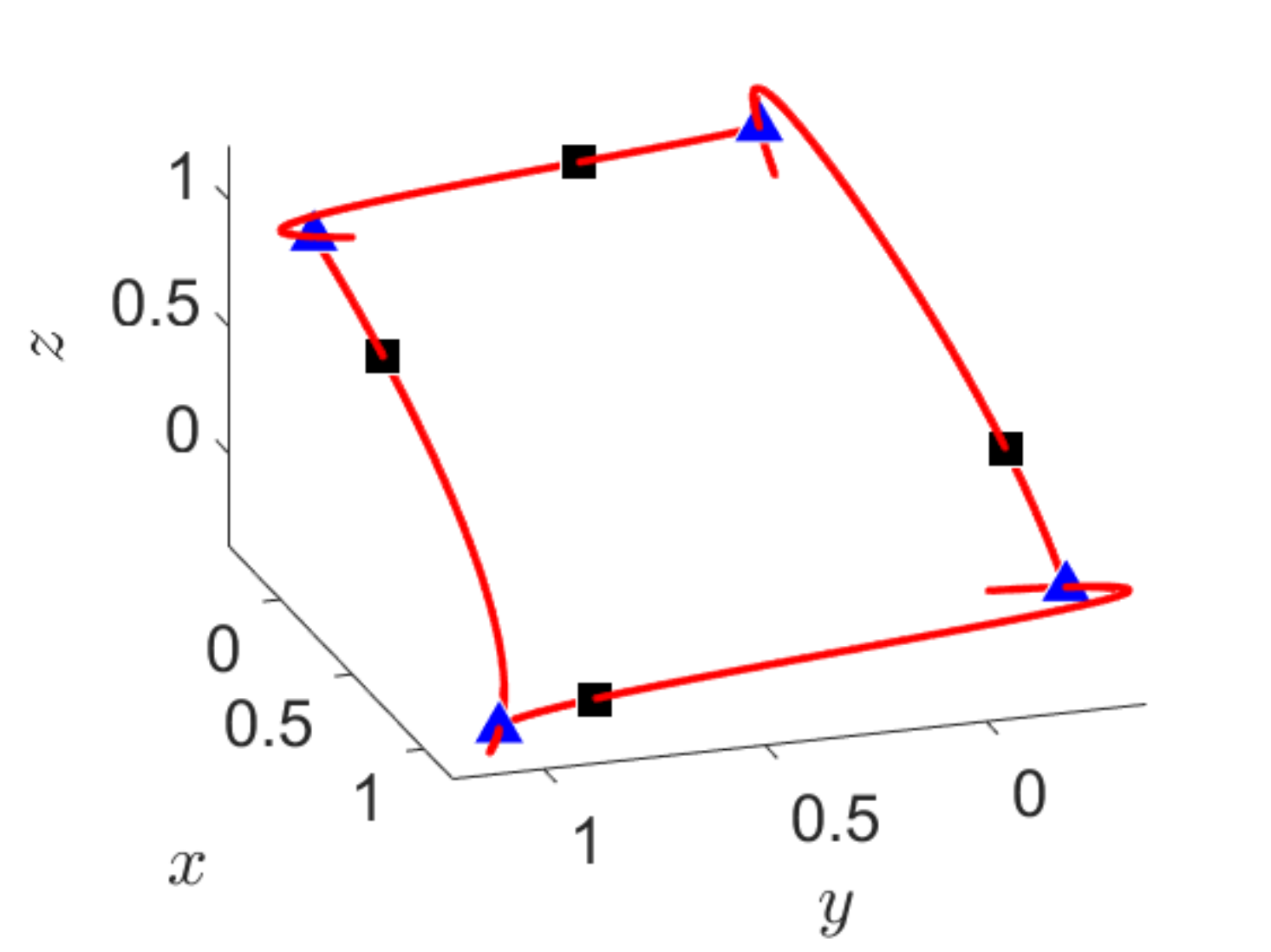}
\caption{Saddle-node connection of stable and saddle mode locked period-4 orbit for $a=1.2, b=0.1$. \label{1loop}}
\end{figure}

Fig. \ref{1loop} shows the saddle period-4 orbit with black squares and the stable period-4 orbit with blue triangles. We observe that they indeed form a saddle-node connection forming a single loop. Thus, the orbit existing for $a=1.2$ is a mode-locked period-4 cycle.

\begin{figure}[!htbp]
\centering
\includegraphics[width=0.49\columnwidth]{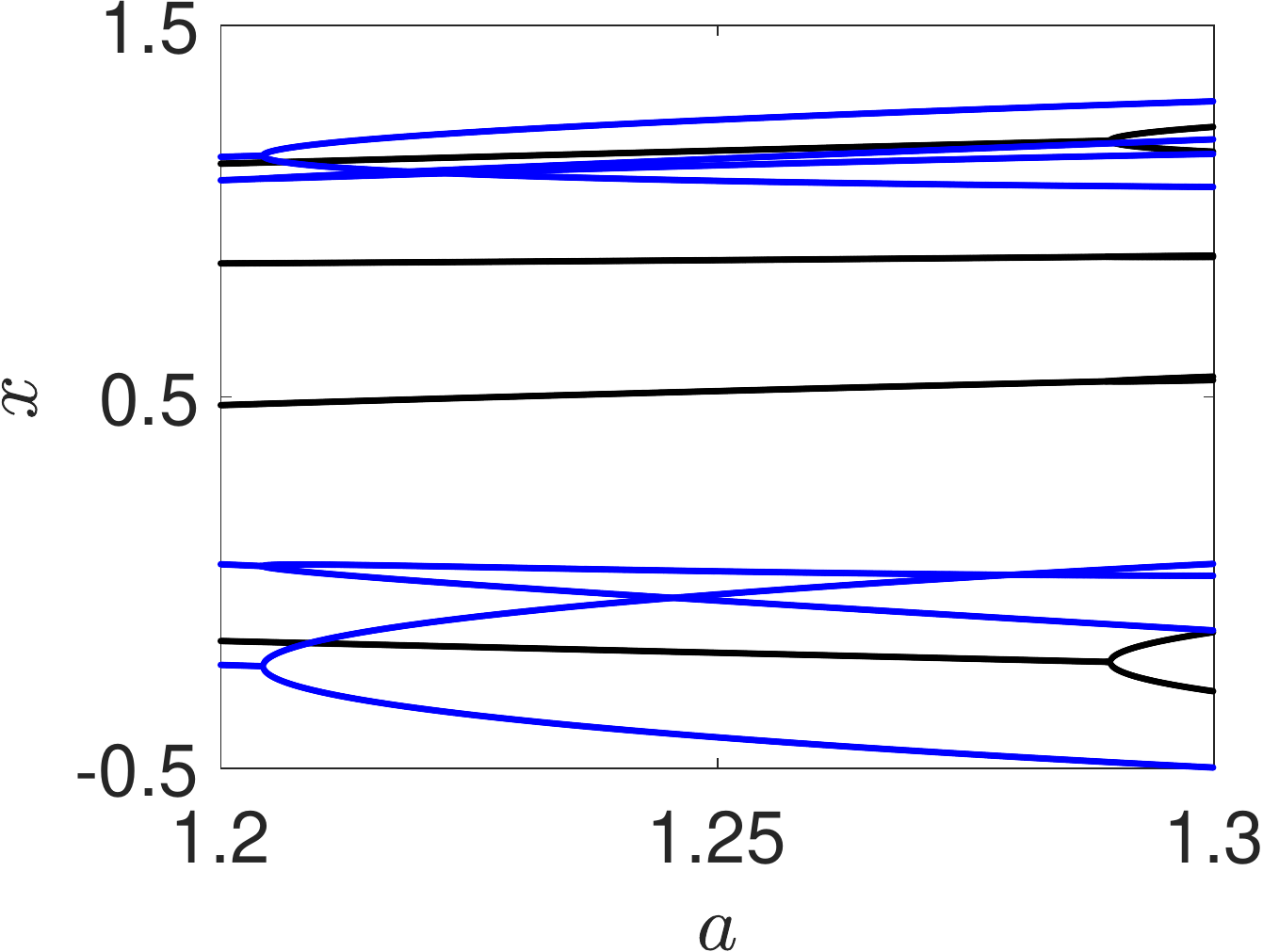}
\caption{Variation of the saddle (black) and stable (blue) periodic orbits with parameter $a$, where $b$ is fixed at 0.1.}
\label{fig:VarSaddlestable_a}
\end{figure}

We now investigate what happens to this cycle as the parameter $a$ is varied. 
A one-parameter bifurcation diagram obtained using a continuation algorithm for both the saddle and stable cycles, is shown in Fig. \ref{fig:VarSaddlestable_a}. It shows that the two orbits bifurcate at different parameter values. There is a range, approximately [1.205, 1.29], where the stable orbit has bifurcated but the saddle has not.

\begin{figure}[!htbp]
\centering
\setlength{\unitlength}{1cm}
\includegraphics[width=0.49\columnwidth]{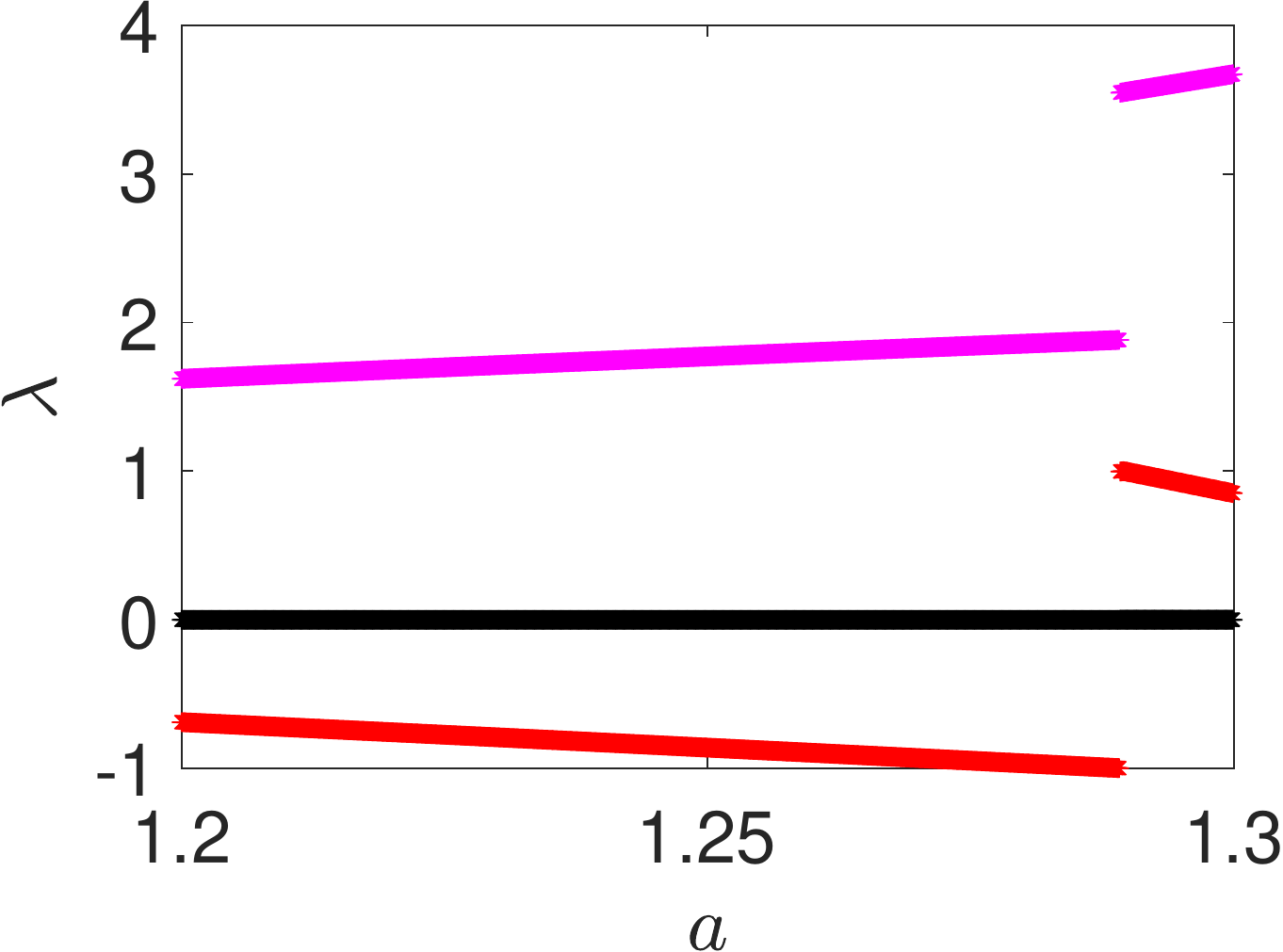}\hfill
\includegraphics[width=0.49\columnwidth]{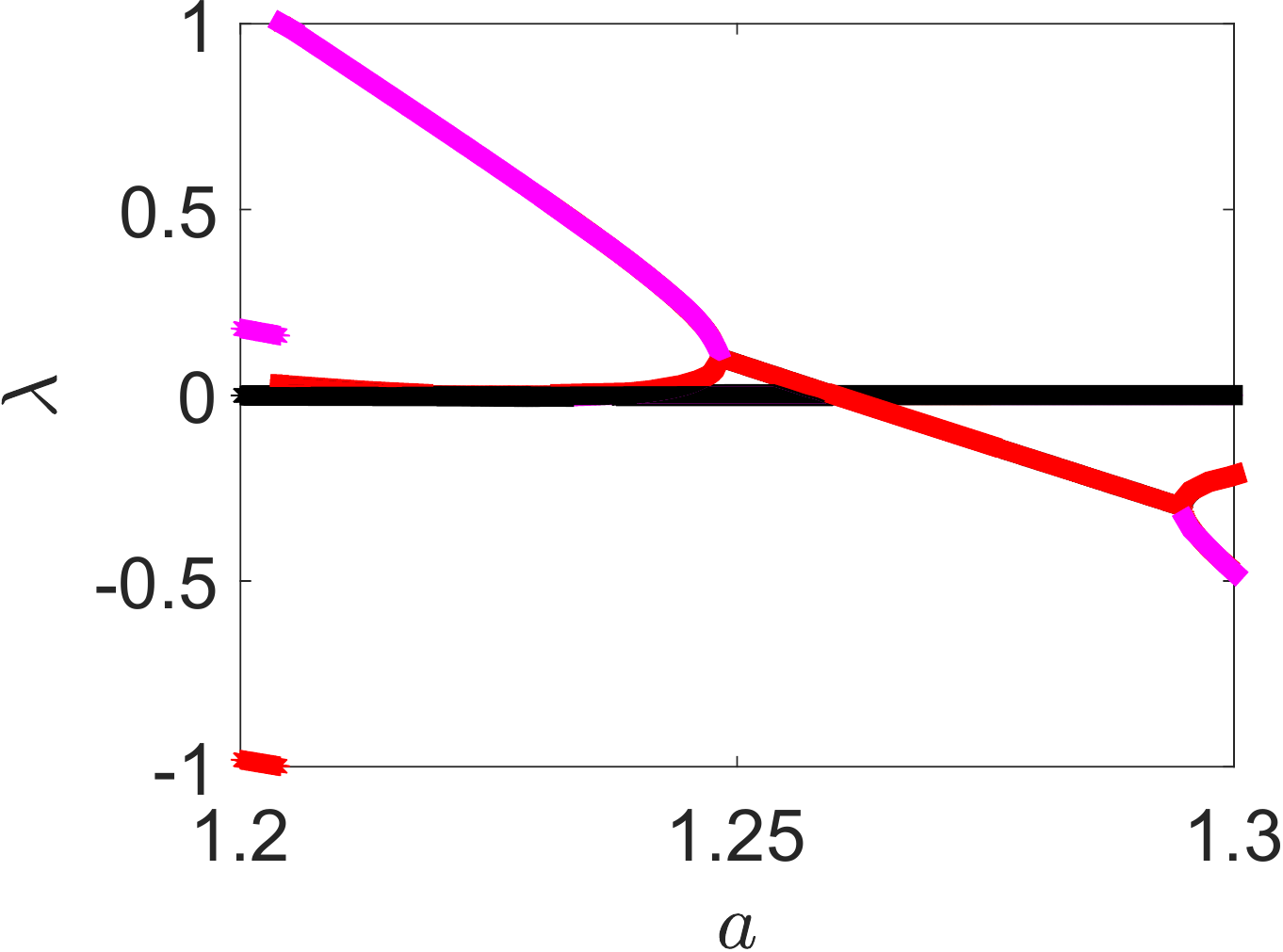}\\
{\footnotesize  (a)} \hspace{1.6in}
{\footnotesize  (b)}
\caption{Variation of the eigenvalue of the saddle cycle (a) and the
  stable cycle (b) with parameter $a$.  The
  parameter $b$ is fixed at 0.1. The three eigenvalues are shown with different colours: $\lambda_1$ is denoted by magenta,
    $\lambda_2$ by red, and $\lambda_3$ by black. }
\label{fig:VarEigenSaddleWith_a}
\end{figure}

Continuation of the eigenvalues of both the saddle and stable
mode-locked periodic orbits as a parameter varies, is presented in
Fig. \ref{fig:VarEigenSaddleWith_a}.
Fig. \ref{fig:VarEigenSaddleWith_a}(a) and (b) show that with increase
of the parameter, at $a=1.204$,  $\lambda_2$ reaches $-1$, which leads to the
period doubling.  

\begin{figure}[!htbp]
\centering
\includegraphics[width=0.47\columnwidth]{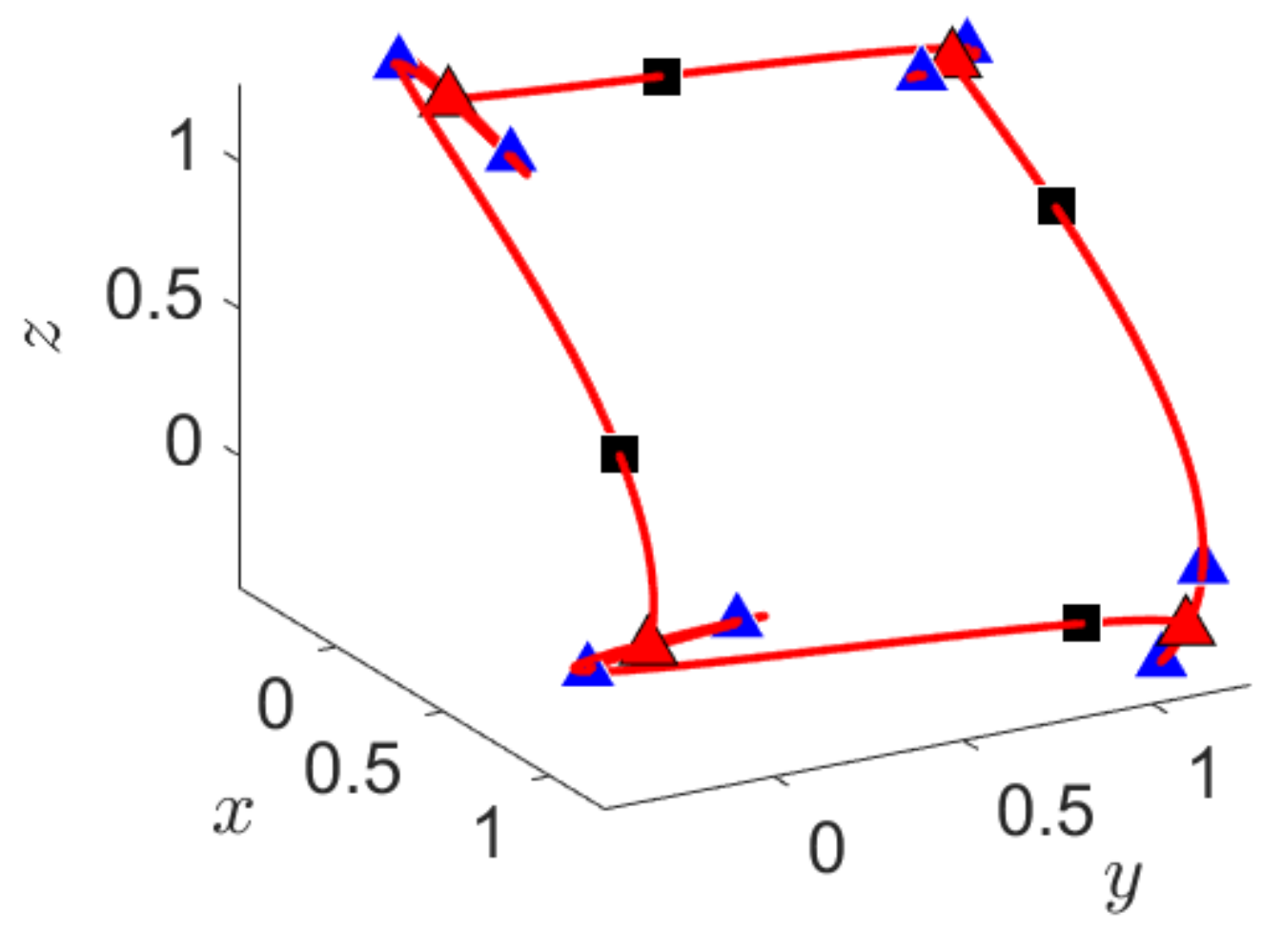}
\includegraphics[width=0.51\columnwidth]{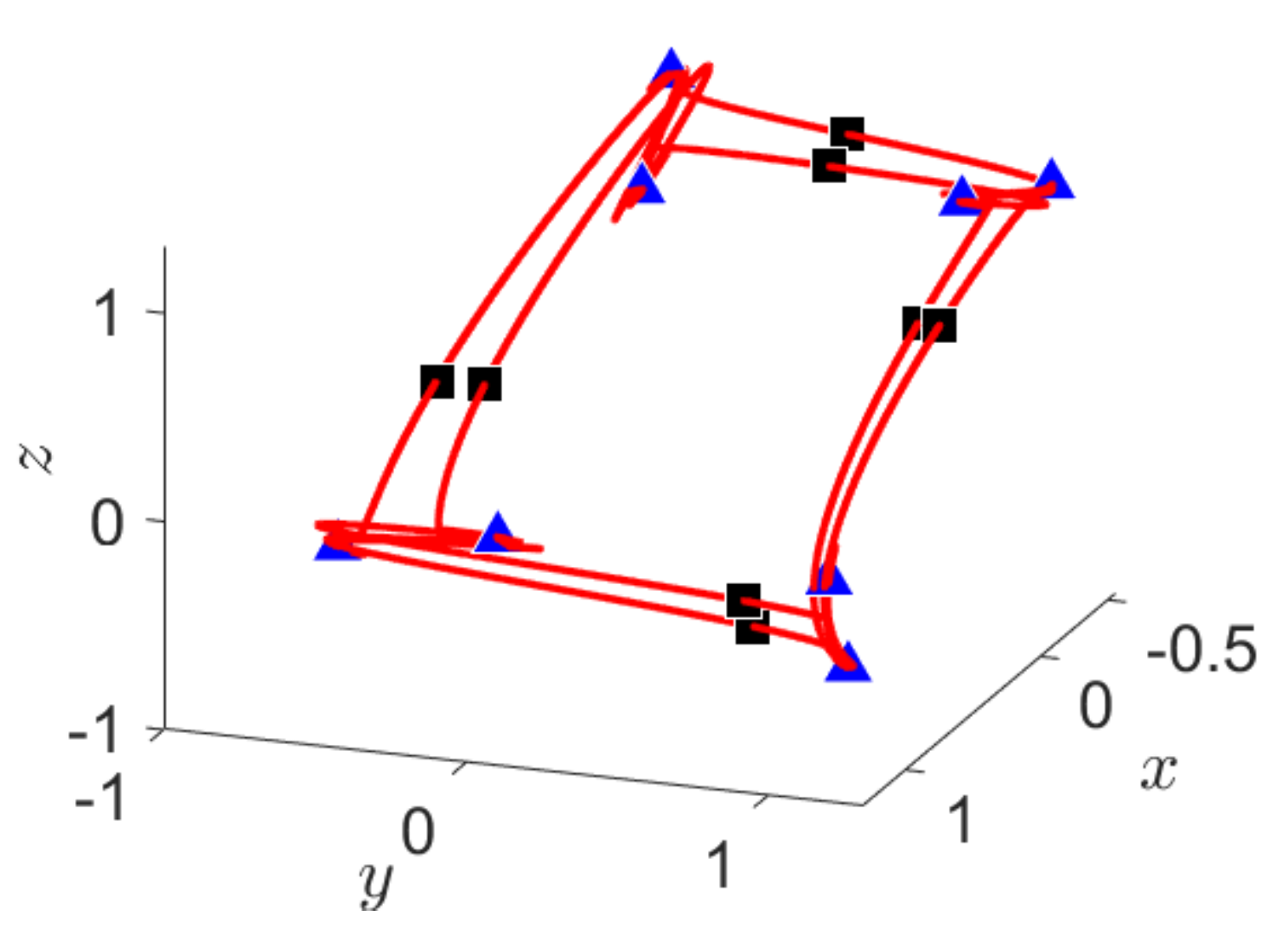}
\\
{\footnotesize  (a)} \hspace{1.6in}
{\footnotesize  (b)}
\caption{(a) The situation at $a=1.25$ and $b=0.1$ where the stable
  period-4 cycle has already doubled but the saddle cycle has not. The
  stable period-8 cycle (shown as blue triangles) has emerged through a flip
  bifurcation of the period-4 cycle (marked in red triangles, now a
  saddle), and the saddle periodic orbit marked by black squares. (b)
  Saddle-node connection of stable and saddle mode locked period-eight
  orbit for $a=1.3, b=0.1$. The red curve denotes the unstable
  manifold making a loop with double length.}
\label{fig:SingleLoopGenHen}

\end{figure}

When the stable cycle has already doubled but the saddle one has not (Fig.~\ref{fig:SingleLoopGenHen}(a)), the unstable manifolds of the saddle points (red line) form a single closed loop structure connecting all the points. Branches emerge from the saddle period-4 cycle (which was earlier the node) to connect the stable period-8 points.

We next consider a parameter value $a=1.3$, $b=0.1$, at which the
saddle and stable period-eight orbits coexist. The period-eight orbit
cyclically visits each of its points. Considering both branches of the
unstable manifold, we observe that the closed invariant curve winds
around twice and its length has doubled, see
Fig. \ref{fig:SingleLoopGenHen}(b).


\begin{table}[tbh]
\begin{center}
\begin{tabular}{ |c|c|c|c| } 
 \hline
 Type & $\lambda_{1}$ & $\lambda_{2}$ & $\lambda_{3}$\\ 
 \hline
 Stable  & 0.1795 & -0.9813 & -0.0006\\
 \hline
 Saddle  & 1.6217 & -0.6890 & -0.0001\\ 
 \hline
\end{tabular}
\end{center}
\caption{The eigenvalues of stable and saddle cycles for $a=1.2, b=0.1$, i.e., before the bifurcation\label{eigens}}
\end{table}


For $a=1.2, b=0.1$, the eigenvalues of the stable and saddle
period-four orbits are given in Table~\ref{eigens}.  It shows that
before the bifurcation, $\lambda_{1}^{N}$ and $\lambda_{1}^{S}$ are
positive and $\lambda_{2}^{N}$ and $\lambda_{2}^{S}$ are negative. 
As the parameter is varied, $\lambda_{2}^{N}$ reaches $-1$ first, followed by
$\lambda_{2}^{S}$. The third eigenvalue is negative. Under this
condition, \cite{Gardini12} conjectured that the shape of the manifold
should be a M\"{o}bius strip and a length-doubling biurcation would
occur. Our numerical results support the conjecture.

\section{From mode-locked orbits to cyclic invariant closed curves}
\label{sec:TransitCyclicClosedCurves}

\begin{figure}[tbh]
  \centering
  \includegraphics[width=2in]{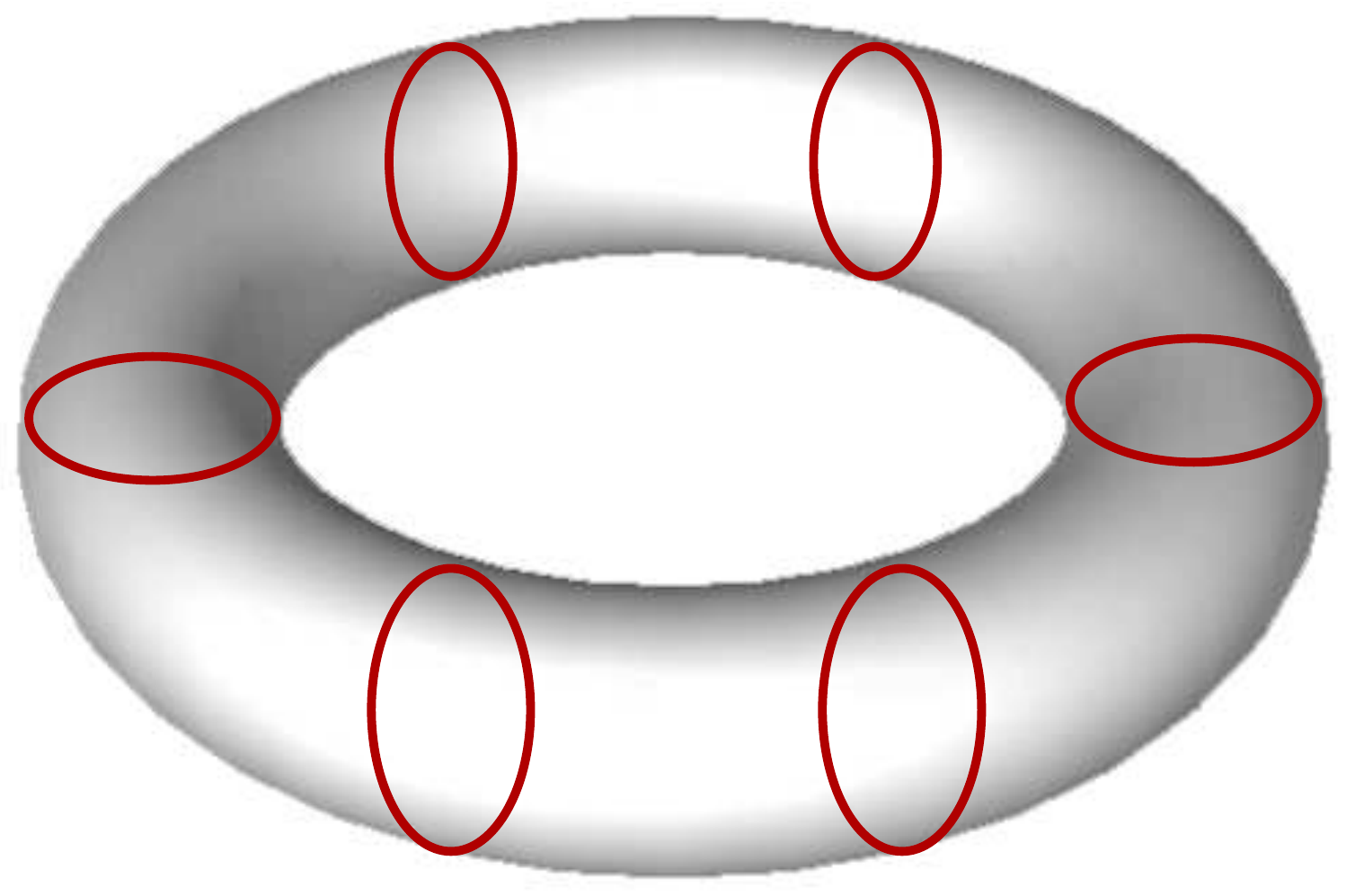}
  \caption{Multiple loops located on a torus}
  \label{fig:torus}
\end{figure}

It is known that when a third frequency is born from a quasiperiodic
orbit, it creates a torus in discrete time.  In this section we
present the mechanisms of generation of a third frequency from a
mode-locked periodic orbit. We find that it can result in two types of
orbits.  The resulting orbit can be a collection of quasiperiodic
loops located on a torus in the discrete-time phase space
(schematically shown in Fig.~\ref{fig:torus}), and the iterates visit
them cyclically. The resulting orbit can also be a mode-locked
periodic orbit where the points are connected by cyclic closed
invariant curves.

\begin{figure*}[tbh]
  \centering
  \includegraphics[width=4.5in,height=1.27in]{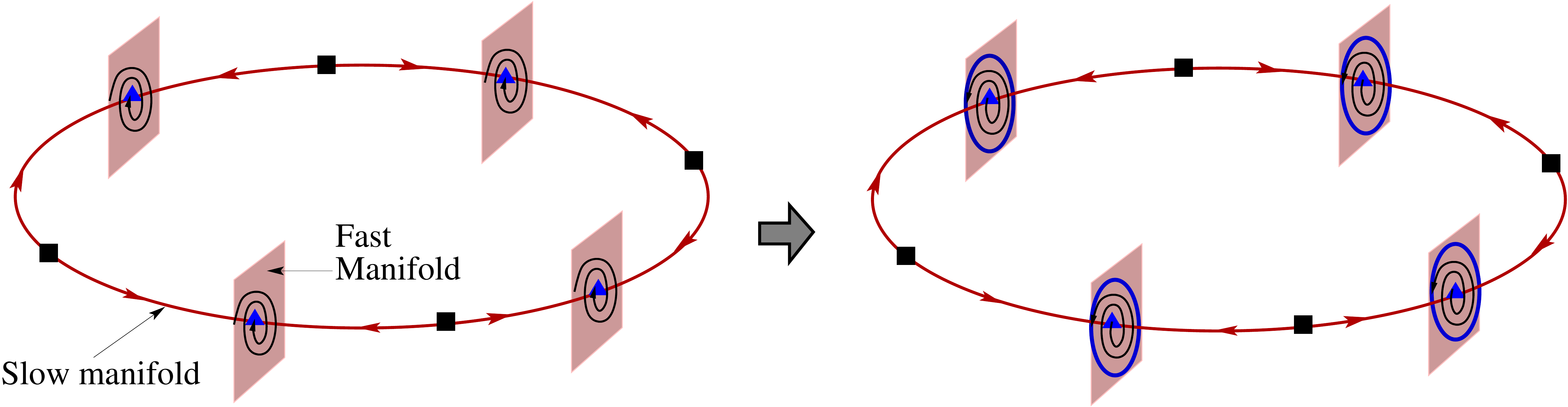}\\{\small (a)}\hspace{2.5in}{\small (b)}\\ \vspace{0.2in}
  \includegraphics[width=4.5in,height=1.27in]{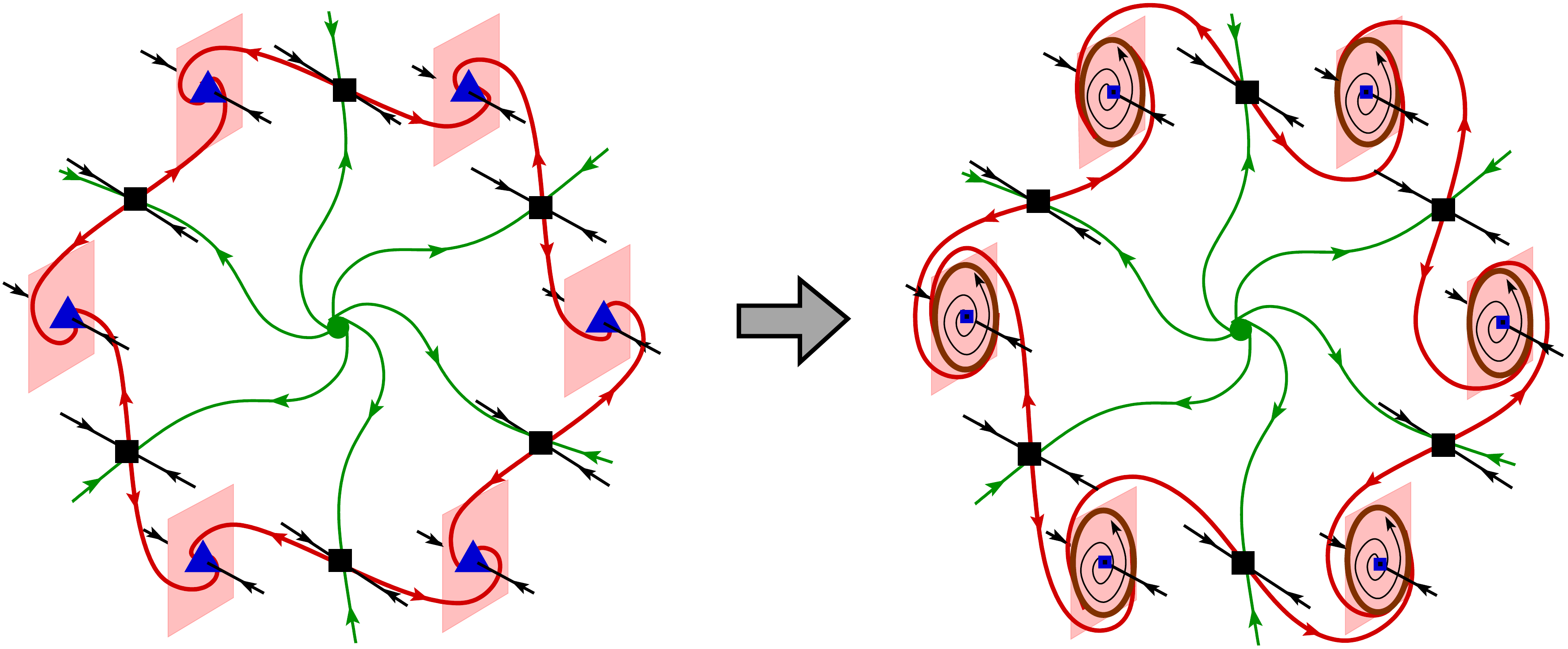}\\{\small (c)} \hspace{2.5in}{\small (d)}
  \caption{Schematic representations of the two mechanisms of creation of multiple loops through Neimark-Sacker bifurcation. The blue triangles represent the foci, the black squares represent the saddles.}
  \label{neimarks}
\end{figure*}

How can such structures be created? Since the birth of a third frequency is caused by a Neimark-Sacker bifurcation, i.e., a pair of complex conjugate eigenvalues exiting the unit circle, the starting point in such a transition has to be a closed invariant curve created out of a saddle-focus connection. If a closed invariant curve is composed of a saddle-node connection as shown in Fig.~\ref{invariant}, with the change of a parameter a pair of real eigenvalues have to turn complex conjugate, thus creating a saddle-focus connection.

The saddle-focus connection can be of two types. The unstable manifold
of the saddle cycle can either connect with the one-dimensional stable
manifold associated with the real eigenvalue or with the
two-dimensional stable manifold associated with the complex conjugate
eigenvalues of the focus. It depends on the rate of convergence along
the manifolds. If the 1D manifold associated with the real eigenvalue
is the slow one, the unstable manifold connects with it (Fig.~\ref{neimarks}(a)). An example of this situation was presented in \cite{de2011local}. If the 2D
manifold is the slow one, the unstable manifold of the saddle spirals
into the node along this manifold (Fig.~\ref{neimarks}(c)).

The Neimark-Sacker bifurcation occurring in these two types of saddle-node connection are shown in Fig.~\ref{neimarks}. In both cases, if there are $n$ points in the stable cycle, $n$ cyclic closed invariant curves are created.

In the first case, the repelling focus creates a loop around it and the saddle-focus connection through the 1D manifold remains intact (see  Fig.~\ref{neimarks}(b)). In the second case, the unstable manifolds of the saddle connect with the closed loops (see  Fig.~\ref{neimarks}(d)).
Note that in bifurcations of saddle-focus loops also, the saddle and the focus may bifurcate at different parameter values. Fig.~\ref{neimarks} shows situations where the focus has bifurcated but the saddle has not. 

We now give a few examples of the birth of cyclic closed invariant curves from saddle-focus connections in physical system models.

\subsection{3D Lotka-Volterra model}

In this section, we provide an explicit example of the situation schematically depicted in Fig.~\ref{neimarks}(c) and (d) using the three-dimensional
Lotka-Volterra map \cite{Gardini87} given by
\begin{equation}
\begin{aligned}
    x' &= x + Rx(1-x-\alpha y - \beta z),\\
    y' &= y + Ry(1-\beta x - y - \alpha z),\\
    z' &= z + Rz(1- \alpha x - \beta y - z),
    \end{aligned}
    \end{equation}
where $R,\alpha, \beta$ are the parameters of the system.

\begin{figure}[b]
\begin{center}
\includegraphics[width=0.5\textwidth]{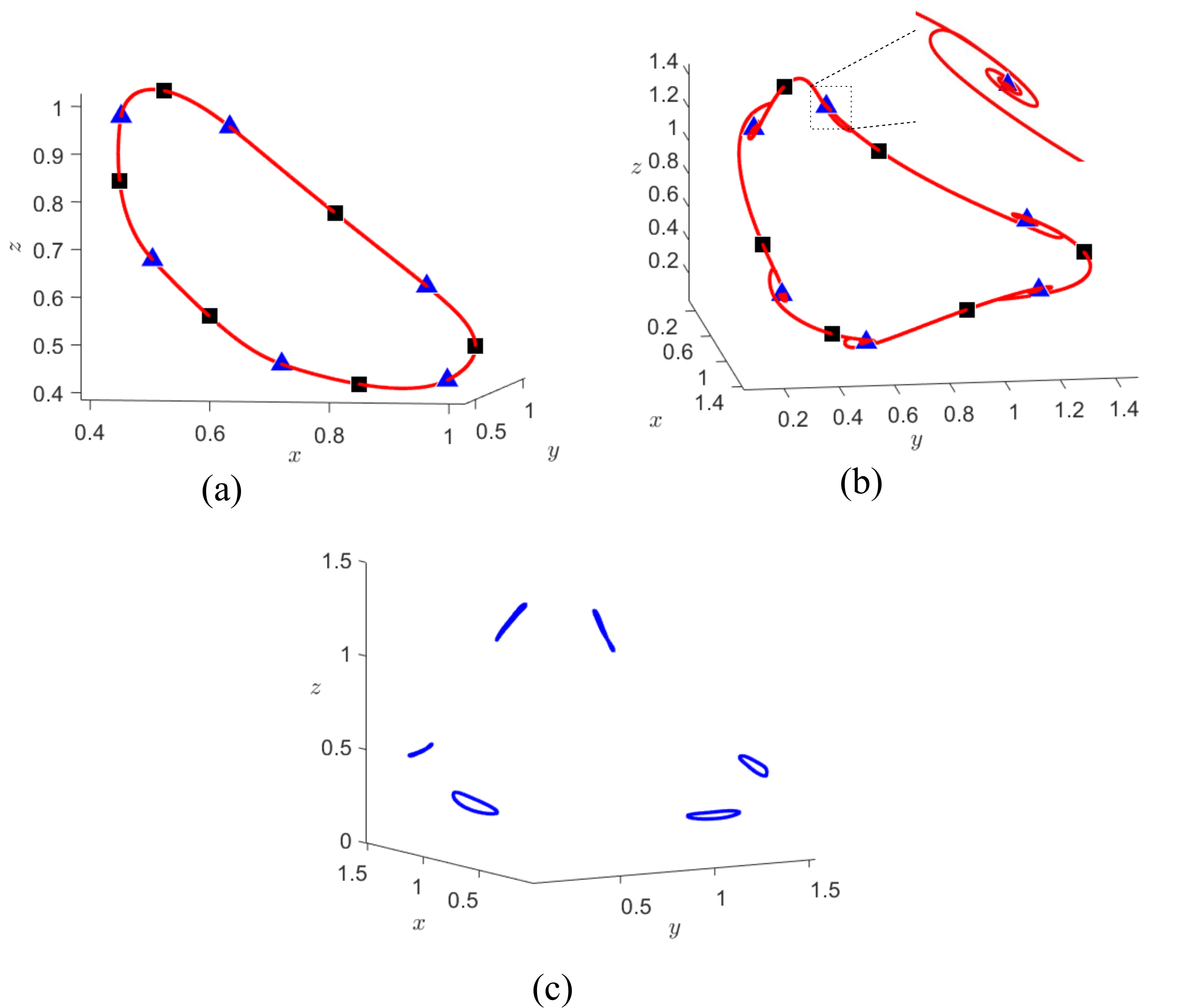}
\caption{In (a), a saddle-node connection
  constituted by a saddle period-six orbit (marked by black squares)
  and a stable period-six orbit (marked by blue triangles) at $\beta = -0.59$. In (b), a
  saddle-focus connection constituted by a saddle period-six orbit
  (marked by black squares) and a stable period-six orbit (marked by blue
  triangles) at $\beta = -0.91$. The unstable manifold of each saddle periodic point is shown in red. In (c), six cyclic invariant
  closed curves are seen at $\beta=-1$ after the period-6 stable cycle has undergone a
  Neimark-Sacker bifurcation.  The other parameters are
  $R=1,\alpha=1$.}
\label{fig:ThreeLVTransitionFull}
\end{center}
\end{figure}

For $R=1, \alpha = 1, \beta = -0.59$, we observe a stable period-six
orbit and a saddle period-six orbit. Considering the one-dimensional
unstable manifold of each saddle periodic point, we observe a
saddle-node connection in Fig. \ref{fig:ThreeLVTransitionFull} (a). As
$\beta$ is decreased to $-0.91$, two eigenvalues of the stable
periodic point become complex conjugate, while the third one remains
real, all of which have modulus less than one. The saddle cycle has
two complex conjugate eigenvalues with modulus less than one and a
real eigenvalue greater than one.  Computation of the one-dimensional
unstable manifolds of each saddle periodic point reveals formation of
a saddle-focus connection (see
Fig. \ref{fig:ThreeLVTransitionFull}(b)).

As the parameter $\beta$ is
further decreased to $-1$, we observe that the period-6 stable cycle has undergone a Neimark-Sacker bifurcation and we can observe
six cyclic closed invariant curves, see
Fig. \ref{fig:ThreeLVTransitionFull}(c).

\begin{figure}[!htbp]
\begin{center}
\includegraphics[width=0.4\textwidth]{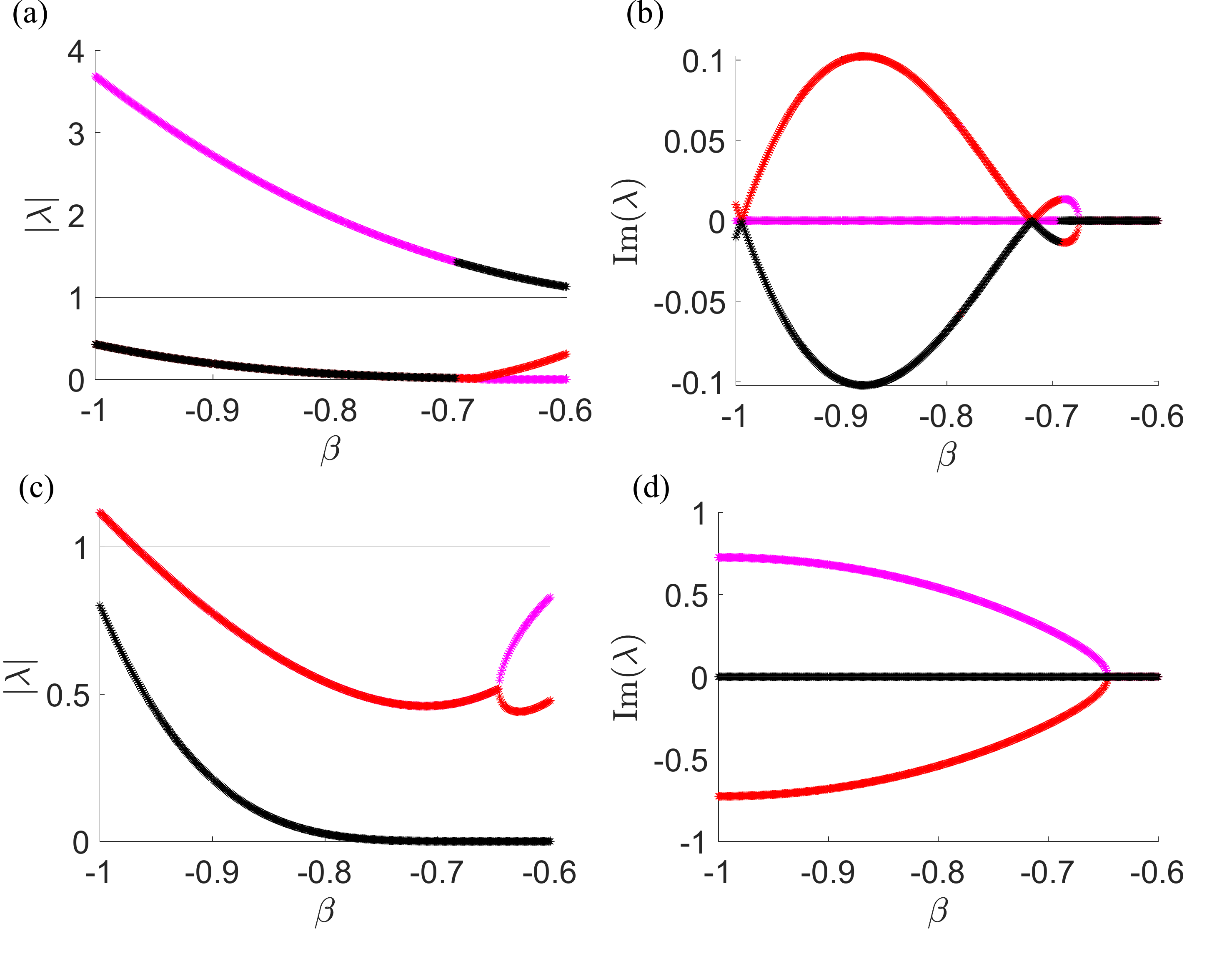}
\caption{The variation of the eigenvalues with respect to the
  parameter $\beta$. (a) the modulus of the eigenvalues of the saddle,
  (b) the imaginary part of the eigenvalues of the saddle, (c) the modulus of the eigenvalues of the node, and (d) imaginary part of the eigenvalues of the node. The three eigenvalues are differentiated by three different colours (magenta, red, black). The other parameters are $R=1,\alpha=1$.}
\label{fig:ThreeLVEigenVariation}
\end{center}
\end{figure}

The variation of the eigenvalues of both saddle and stable periodic
point with respect to the parameter $\beta$ is shown in
Fig. \ref{fig:ThreeLVEigenVariation}.  It shows that two eigenvalues of the stable cycle become complex at $\beta = -0.646$. The modulus crosses $1$ at $\beta=-0.968$. Subsequently we see the onset of a six cyclic closed invariant curve.

\subsection{Three-dimensional border collision normal form}

Piecewise smooth maps occur in many physical and engineering
systems. It has been shown that, in the neighborhood of a
border-crossing fixed point, such systems are aptly represented by a
piecewise linear `normal form' map \cite{BCNF99}. For three-dimensional piecewise
smooth systems, the 3D border collision normal form is given by 
\begin{equation}
X' = F(X) = 
     \begin{cases} 
      A_{L}X + b \mu ,& x\leq 0 \\
      A_{R}X + b \mu ,& x > 0, \\
        \end{cases}
\end{equation}
where $F : \mathbb{R}^{3} \rightarrow \mathbb{R}^{3}, X = [x,y,z]^{T},b = [1,0,0]^{T}$ is a column vector, and $\mu \in \mathbb{R}$. The matrices $A_{L}, A_{R}$ are defined by 
\begin{equation}
A_{L} = 
\begin{bmatrix}
\tau_{L} & 1 & 0\\
-\sigma_{L} & 0 & 1\\
\delta_{L} & 0 & 0
\end{bmatrix},\;\;\;
\mbox{and} \;\;\;
A_{R} = 
\begin{bmatrix}
\tau_{R} & 1 & 0\\
-\sigma_{R} & 0 & 1\\
\delta_{R} & 0 & 0
\end{bmatrix}.
\end{equation}

\begin{figure}[!htbp]
\begin{center}
\includegraphics[width=0.45\textwidth]{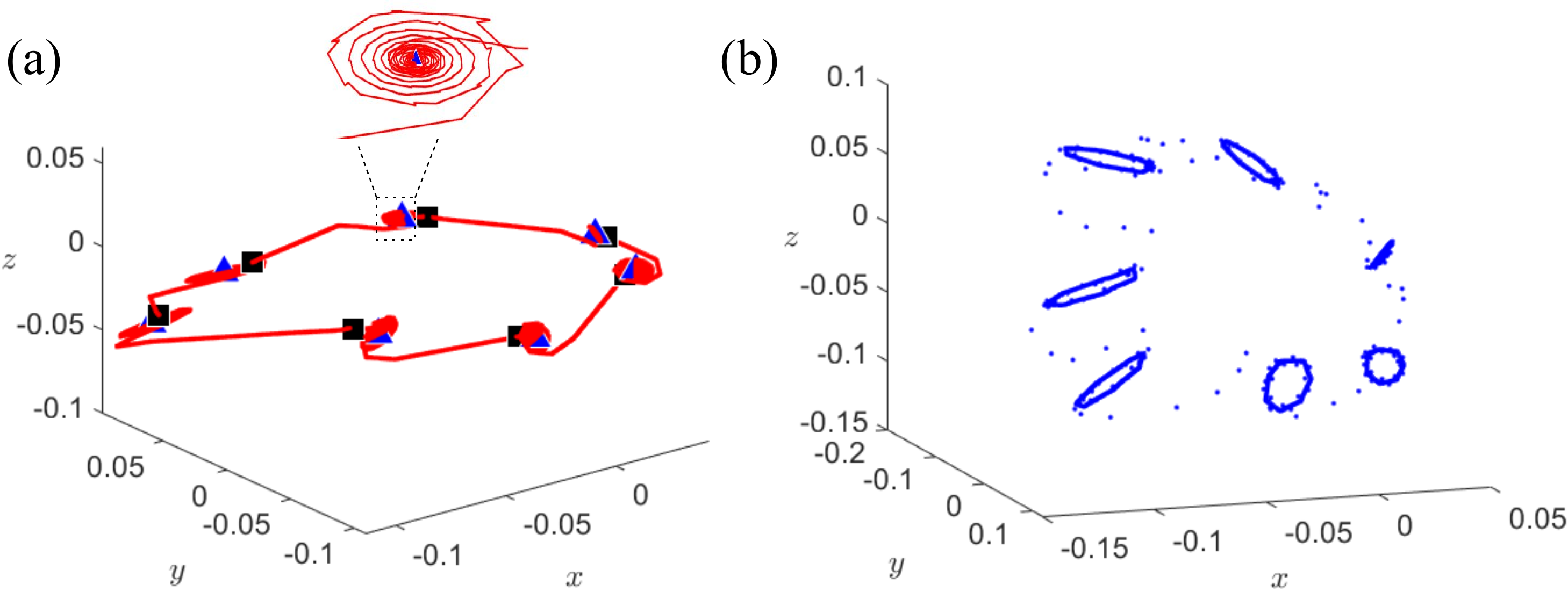}
\caption{In (a), a saddle-focus connection at $\delta_{R}=1.4$ constituted by period-seven saddle orbit (marked by black squares) and stable period-seven orbit (marked by blue triangles). The unstable manifold of each saddle periodic point is shown in red.  In (b), seven cyclic invariant closed curves are seen at $\delta_{R}=1.5$ after each stable periodic point has undergone a Neimark-Sacker bifurcation. The other parameters are $\tau_{R}=-0.5, \tau_{L}=0.74, \delta_{L}=0.73, \mu=0.01, \sigma_{L} = 0.5, \sigma_{R} = 1.1$. }
\label{fig:ThreeBCNFSadlFocusLoop}
\end{center}
\end{figure}

\begin{figure}[!htbp]
\begin{center}
\includegraphics[width=0.45\textwidth]{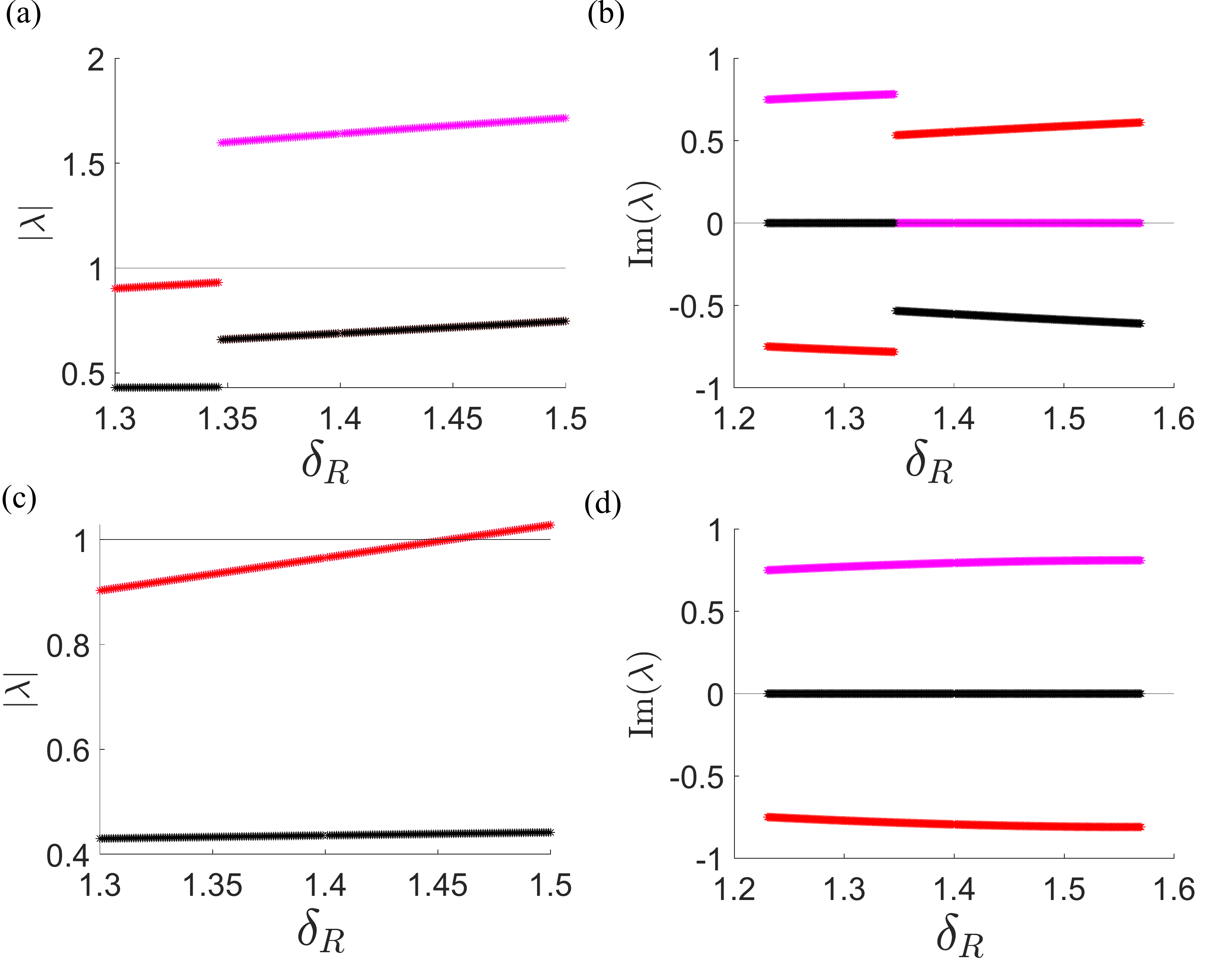}
\caption{The variation of the eigenvalues with respect to the parameter $\delta_{R}$.
In (a), the modulus of the eigenvalues of the saddle, (b) the imaginary part of the eigenvalues of the saddle,
(c) the modulus of the eigenvalues of the node, and (d) imaginary part of the eigenvalues of the node. 
The three eigenvalues are differentiated by three different colours (magenta, red, black). The other parameters are set as $\tau_{R} = -0.5, \tau_{L} = 0.74, \delta_{L} = 0.73, \mu = 0.01, \sigma_{L} = 0.5, \sigma_{R}= 1.1$.}
\label{fig:BCNFTransitionFull}
\end{center}
\end{figure}

In Fig.\ref{fig:ThreeBCNFSadlFocusLoop} (a), we observe a period-7
stable and saddle orbit marked by blue triangles and black squares
respectively for $\delta_{R} = 1.4$. A saddle-focus connection is
observed. After increasing $\delta_{R}$ to $1.5$, seven cyclic closed
invariant curves are observed. Continuation of the eigenvalues with
variation of the parameter $\delta_{R}$
(Fig. \ref{fig:BCNFTransitionFull}) shows that for
$\delta_{R}=1.4$, both the saddle cycle and the stable cycle have two
complex conjugate eigenvalues and a real eigenvalue. At $\delta_{R}=1.455$, the complex conjugate eigenvalues cross the unit circle, and a 7-piece cyclic closed invariant curve is born.

\section{Birth of three-frequency resonant torus}
\begin{figure}
\centering
\includegraphics[width=0.5\textwidth]{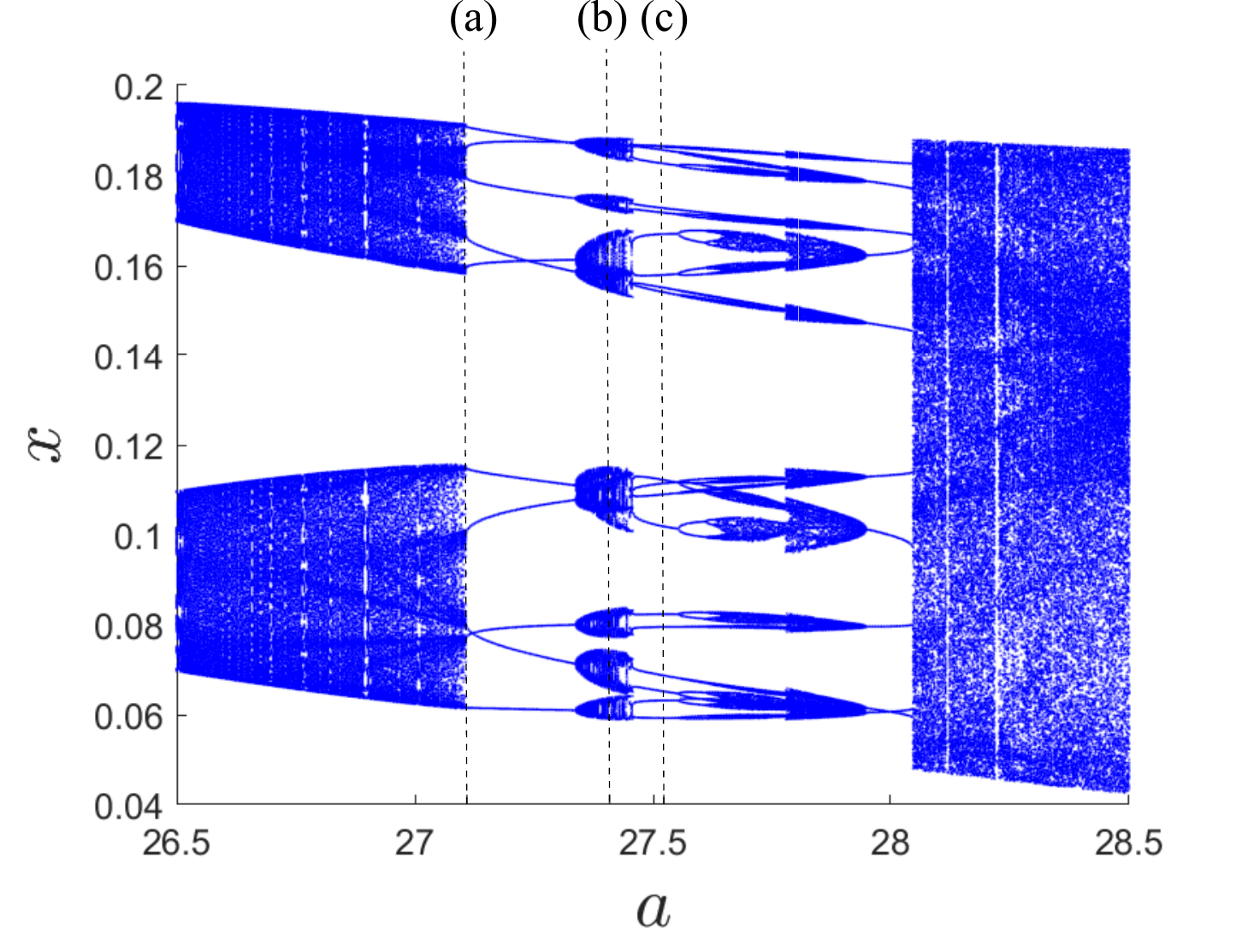}
\caption{ A one-parameter bifurcation diagram of \eqref{eq:globalcoupled} with respect to parameter $a$. The other parameters
are fixed as
$p=0.5, \epsilon = -1.4, b= 1.688, c=3.5, d=0.85$. Three different parameter values are chosen as shown by dotted lines. }
\label{fig:BifurThreeFreq}
\end{figure}

Let us consider the globally coupled three-dimensional map \eqref{eq:globalcoupled}. A similar two-dimensional map has been studied in \cite{Zhu08}. The dynamical equations of the map are given by
\begin{equation}
\begin{aligned}
    x' &= f(x) + p \epsilon \Bigl(f(y) - f(x) \Bigr),\\
    y' &= f(y) + (1-p) \epsilon \Bigl(f(x) - f(y)\Bigr),\\
    z' &= x,
    \end{aligned}
    \label{eq:globalcoupled}
    \end{equation}
    where
    \begin{equation*}
      f (x) = x(1-x)(ax^{2} + (b^2 - da)x + c),
    \end{equation*}
and $ p, \epsilon, a, b, c, d$ are parameters. A one-parameter
bifurcation diagram with respect to parameter $a$ is presented in Fig~\ref{fig:BifurThreeFreq}. 
For $a < 27.1107$, two disjoint cyclic quasiperiodic closed invariant
curves exist.  At $a = 27.1107$, a period-$10$ orbit is born via a
saddle-node bifurcation. Fig. \ref{fig:CollageThreeFreq}(a) shows that
the period-10 cycle occurs on two disjoint closed invariant curves
connected by the one-dimensional unstable manifolds of a period-10 saddle
cycle.  When $a = 27.408$, we observe each of the stable
periodic points undergo a supercritical Neimark-Sacker bifurcation
resulting in the formation of ten-disjoint cyclic closed loops, see
Fig. \ref{fig:CollageThreeFreq}(b).  

 The system develops two coexisting stable period-$20$
  orbit for $a \in (27.461, 27.6)$ (this is not shown to avoid the
  figure becoming messy). At $a = 27.46$, ten disjoint cyclic closed
  loops are formed via a saddle-node bifurcation of a pair of stable
  and saddle period-$20$ orbits. Fig. \ref{fig:SaddleNodeP20} shows
  the bifurcation diagram in this parameter range. For clarity, we
  show the occurrence of saddle-node bifurcation for the cycles lying
  on a single loop (two stable period-2 and two saddle period-2
  cycles). The concerned loop is shown in
  Fig. \ref{fig:CollageThreeFreq}(d).

  Similar saddle-node bifurcations occur in all 10 loops. 
 At $a = 27.5210$, there are 10 closed invariant curves formed through saddle-node connections. These are, in turn, located on two loops. The iterates toggle between the two bigger loops and move cyclically among the smaller loops. The resulting phase space is shown in
Fig. \ref{fig:CollageThreeFreq}(c). Note that the saddle and stable
periodic orbits form cyclic disjoint loops connected via their
one-dimensional unstable manifold. A zoomed-in version of one of the
cyclic closed loops is shown in Fig. \ref{fig:CollageThreeFreq}(d).

\begin{figure}[!htbp]
\centering
\includegraphics[width=0.48\textwidth]{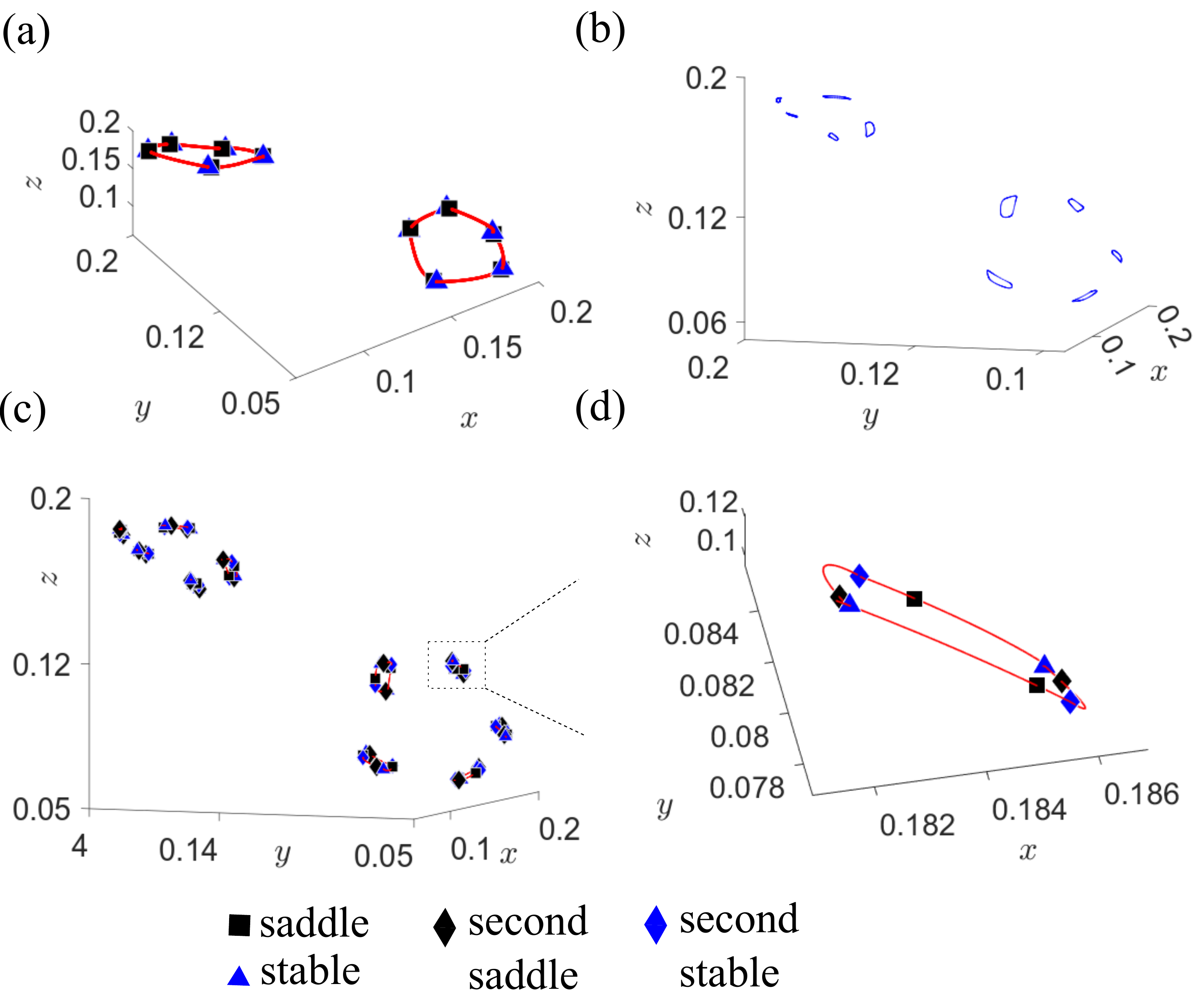}
\caption{In (a), two disjoint cyclic closed invariant loops are
  observed comprising a stable period-10 orbit and a saddle period-10
  orbit at $a = 27.1107$. In (b), ten disjoint pieces of cyclic closed
  invariant curves are observed at $a = 27.408$. In (c), ten disjoint
  pieces of closed invariant loops are observed comprising two
  stable period-$20$ orbits and two saddle period-$20$ orbits
  at $a=27.5210$. In (d), a magnified version of a single closed loop
  is shown. }
\label{fig:CollageThreeFreq}
\end{figure}

\begin{figure}[!htbp]
\centering
\includegraphics[width=0.3\textwidth]{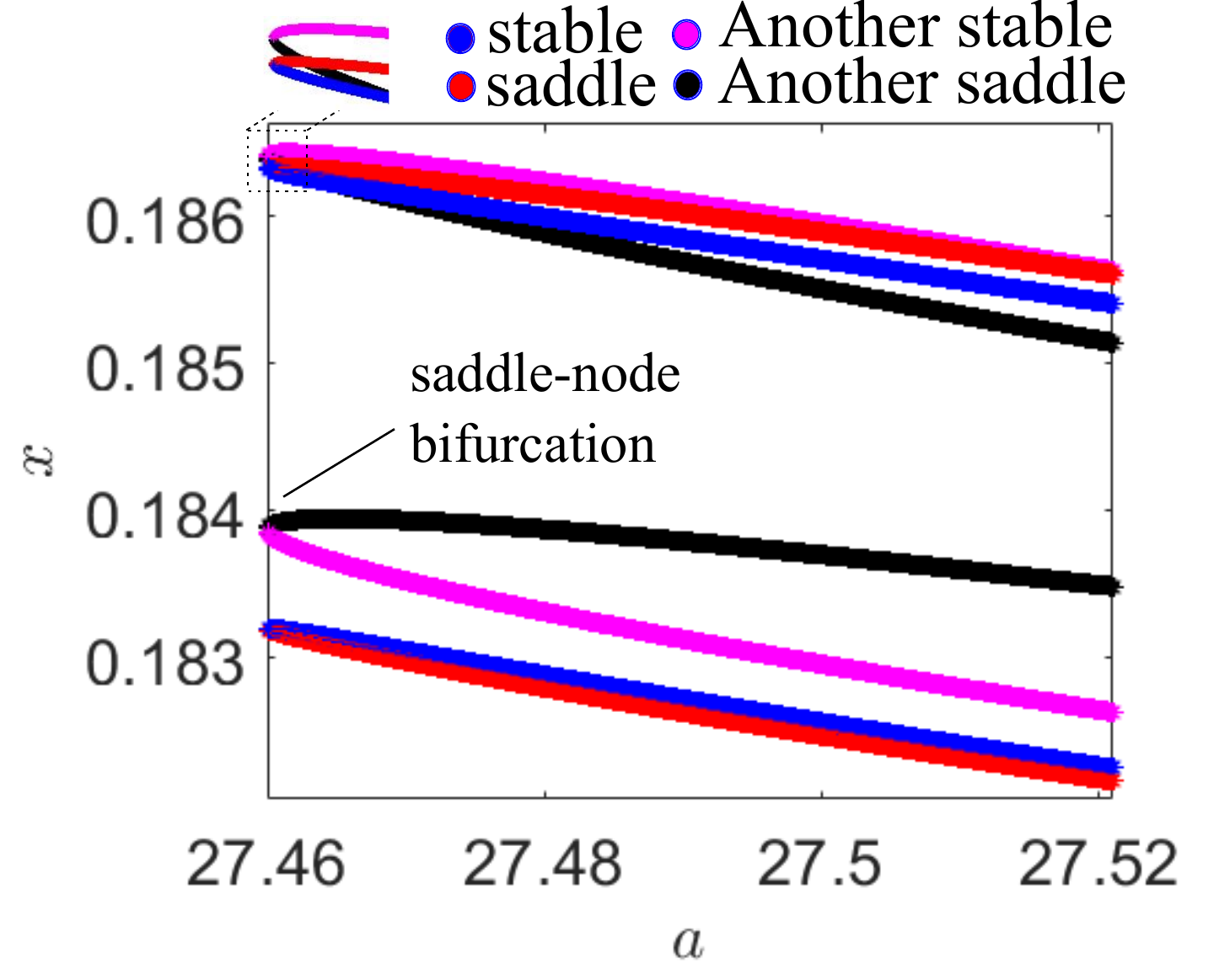}
\caption{Continuation of periodic points lying on a
    single cyclic closed loop with respect to parameter
    $a$. Saddle-node bifurcation of a pair of two stable and two
    saddle periodic points lying on a single loop at $a=27.46$. Two
    stable and saddle points are marked in blue and red
    respectively. The other two coexisting stable and saddle points
    are marked in magenta and black respectively. The other parameters
    are fixed as $p=0.5, \epsilon = -1.4, b= 1.688, c=3.5, d=0.85$.}
\label{fig:SaddleNodeP20}
\end{figure}

\section{Collision of saddle-node and saddle-saddle connections}
\label{sec:collision}

It has been reported in \cite{Banerjee12}, that there can be
situations where a stable torus and an unstable torus collide and
disappear (see Fig.~\ref{quasifig}(d)). If the parameter is varied in
the opposite direction, one can see the birth of a torus out of
nothing. Such bifurcations, involving ergodic tori have been reported
in power electronic systems.

\begin{figure}[!htbp]
\begin{center}
\includegraphics[width=0.3\textwidth]{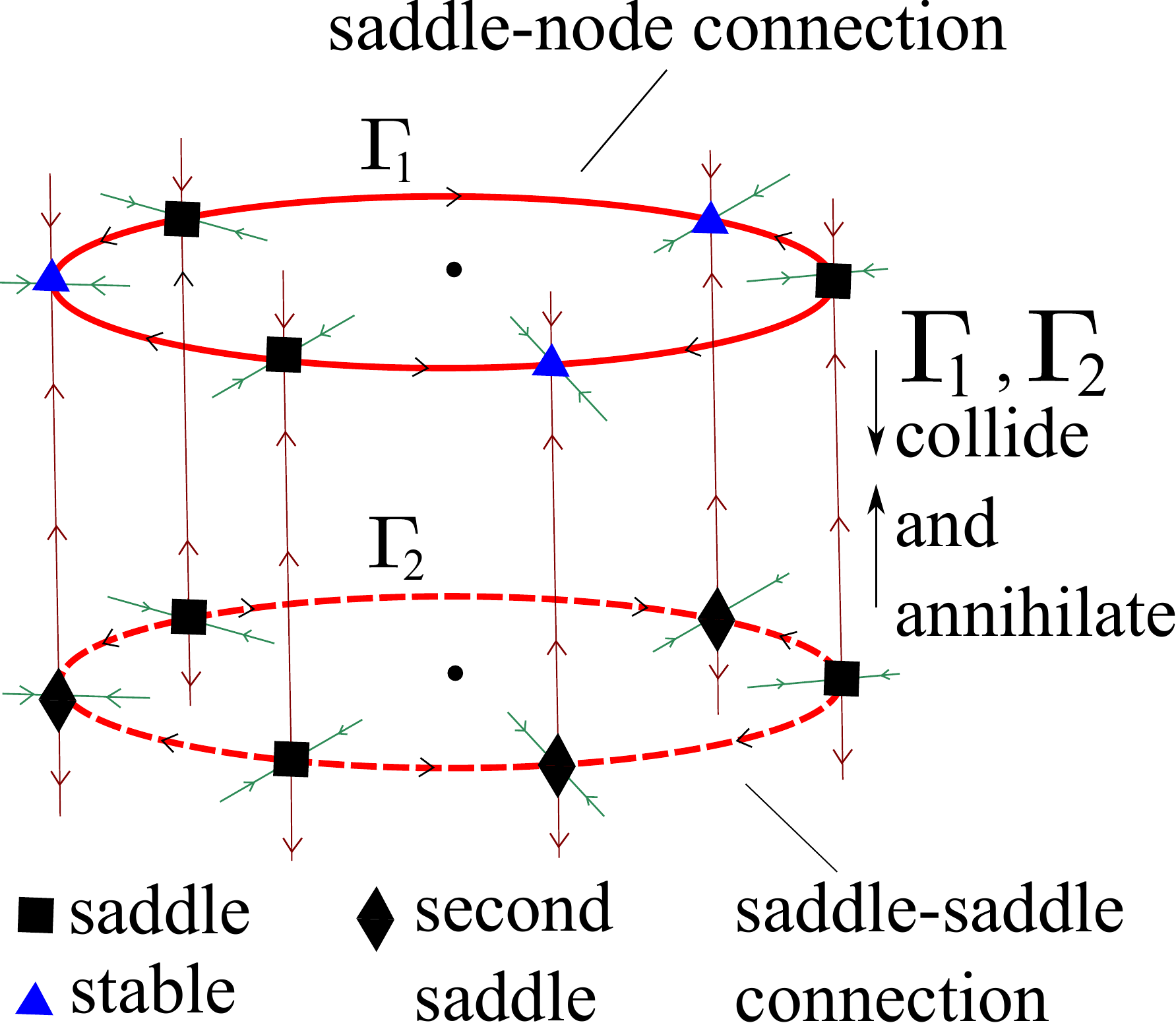}
\caption{Collision of saddle-node connection $\Gamma_{1}$ with the saddle-saddle connection $\Gamma_{2}$. The saddle periodic orbits are denoted by black squares, another saddle periodic orbit with black diamond, and the stable periodic orbits are denoted by blue triangles. }
\label{fig:CollidingConnection}
\end{center}
\end{figure}

The question is, can such a bifurcation involve resonant tori, i.e.,
mode-locked periodic orbits? In Fig.~\ref{fig:CollidingConnection},
we schematically show how such a collision between saddle-node
connection $\Gamma_{1}$ and saddle-saddle connection $\Gamma_{2}$ can
take place. The three eigendirections are shown in different colours
and arrows.

As $\Gamma_{1}$ and $\Gamma_{2}$ approach each other as the parameter
varies, they may collide and annihilate each other. If the parameter
varies in the opposite direction, there is a sudden appearance of two
closed invariant curves, one with a saddle-node connection and the
other with a saddle-saddle connection. This is a topologically
feasible scenario, but we have not yet found any physical system that
exhibits this phenomenon.

\section{Conclusions}
\label{sec:Conclusions}

In this work, we have explored four scenarios related to the bifurcations of mode-locked periodic orbits and their associated closed invariant curves:
\begin{enumerate}
\item The birth of two disjoint closed invariant curves
\item The doubling of the length of the closed invariant curve
\item The birth of a few pieces of cyclic closed invariant curves
\item Merger and disappearance of two closed invariant curves, one involving a saddle-node connection and the other involving a saddle-saddle connection.
\end{enumerate}

We have provided prototypical examples of three-dimensional maps in
which bifurcation phenomena 1, 2, and 3 above can be observed. Through
numerical simulations, we have validated the conjectures proposed in
\cite{Gardini12}. In addition, we have shown the interesting
structures of the invariant manifolds when the node has bifurcated and
the saddle has not, and vice versa.

We have explored the transitions from mode-locked periodic orbits to
the formation of cyclic closed invariant curves.  There can be two
mechanisms of the creation of an $n$-piece closed invariant curve out
of a period-$n$ mode locked orbit. We have presented examples of such
bifurcations from a saddle-focus connection.

The situation depicted in
Fig.~\ref{fig:CollidingConnection} is mathematically possible, where
two closed invariant curves---one with a saddle-node connection and
the other with a saddle-saddle connection---merge and disappear. We
have not yet found physical examples of such bifurcations, which
remains an open problem.

\section*{Acknowledgements}

SSM expresses his thanks to Dr. David J.W. Simpson and Prof. Hil
G. E. Meijer for many insightful suggestions in the course of this
work. SSM also acknowledges the IISER Kolkata post-doctoral fellowship
for financial support. SB acknowledges the J C Bose National
Fellowship provided by SERB, Government of India, Grant
No. JBR/2020/000049.

\section*{Data Availability}
The data that support the findings of this study are available from the corresponding author upon request.
\bibliographystyle{unsrt}
\bibliography{refer}

\begin{thebibliography}{10}

\bibitem{Kaneko83}
K.~Kaneko.
\newblock {Doubling of Torus}.
\newblock {\em Progress of Theoretical Physics}, 69(6):1806--1810, 1983.

\bibitem{Spiegel83}
A.~Arnéodo, P.H. Coullet, and E.A. Spiegel.
\newblock Cascade of period doublings of tori.
\newblock {\em Physics Letters A}, 94(1):1--6, 1983.

\bibitem{Zhu21}
Z.~T. Zhusubaliyev, V.~Avrutin, and A.~Medvedev.
\newblock Doubling of a closed invariant curve in an impulsive {G}oodwin’s
  oscillator with delay.
\newblock {\em Chaos, Solitons \& Fractals}, 153:111571, 2021.

\bibitem{Banerjee12}
S.~Banerjee, D.~Giaouris, P.~Missailidis, and O.~Imrayed.
\newblock Local bifurcations of a quasiperiodic orbit.
\newblock {\em International Journal of Bifurcation and Chaos}, 22(12):1250289,
  2012.

\bibitem{Gardini87}
L.~Gardini, R.~Lupini, C.~Mammana, and M.~G. Messia.
\newblock Bifurcations and transitions to chaos in the three-dimensional
  {L}otka–{V}olterra map.
\newblock {\em SIAM Journal on Applied Mathematics}, 47(3):455--482, 1987.

\bibitem{ding2004interaction}
W-C Ding, JH~Xie, and QG~Sun.
\newblock Interaction of {H}opf and period doubling bifurcations of a
  vibro-impact system.
\newblock {\em Journal of Sound and Vibration}, 275(1-2):27--45, 2004.

\bibitem{Ash95}
P.~Ashwin and J.~W. Swift.
\newblock Torus doubling in four weakly coupled oscillators.
\newblock {\em International Journal of Bifurcation and Chaos},
  05(01):231--241, 1995.

\bibitem{Anish05}
V.~Anishchenko, S.~Nikolaev, and G.~Strelkova.
\newblock Oscillator of quasiperiodic oscillations. two-dimensional torus
  doubling bifurcation.
\newblock {\em International Symposium on Nonlinear Theory and its
  Applications}, pages 23--25, 2005.

\bibitem{Seki01}
M.~Sekikawa, T.~Miyoshi, and N.~Inaba.
\newblock Successive torus doubling.
\newblock {\em IEEE Transactions on Circuits and Systems I: Fundamental Theory
  and Applications}, 48(1):28--34, 2001.

\bibitem{Seki21}
M.~Sekikawa and N.~Inaba.
\newblock Chaos after accumulation of torus doublings.
\newblock {\em International Journal of Bifurcation and Chaos}, 31(01):2150009,
  2021.

\bibitem{Seleznev2}
B.~P. Bezruchko, S.P. Kuznetsov, and Y.P. Seleznev.
\newblock Experimental observation of dynamics near the torus-doubling terminal
  critical point.
\newblock {\em Phys. Rev. E}, 62:7828--7830, 2000.

\bibitem{forcedlogistic}
S.~Kuznetsov, U.~Feudel, and A.~Pikovsky.
\newblock Renormalization group for scaling at the torus-doubling terminal
  point.
\newblock {\em Phys. Rev. E}, 57:1585--1590, 1998.

\bibitem{Meiss06}
D.~B. Wysham and J.D. Meiss.
\newblock Iterative techniques for computing the linearized manifolds of
  quasiperiodic tori.
\newblock {\em Chaos: An Interdisciplinary Journal of Nonlinear Science},
  16(2):023129, 2006.

\bibitem{Vitolo08}
H.~Broer, C.~Simó, and R.~Vitolo.
\newblock Hopf saddle-node bifurcation for fixed points of 3d-diffeomorphisms:
  Analysis of a resonance ‘bubble’.
\newblock {\em Physica D: Nonlinear Phenomena}, 237(13):1773--1799, 2008.

\bibitem{Zhu06}
Zhanybai~T. Zhusubaliyev and Erik Mosekilde.
\newblock Birth of bilayered torus and torus breakdown in a piecewise-smooth
  dynamical system.
\newblock {\em Physics Letters A}, 351(3):167--174, 2006.

\bibitem{Gon21}
A.~S. Gonchenko, S.~V. Gonchenko, and D.~Turaev.
\newblock Doubling of invariant curves and chaos in three-dimensional
  diffeomorphisms.
\newblock {\em Chaos: An Interdisciplinary Journal of Nonlinear Science},
  31(11):113130, 2021.

\bibitem{grebogi1983three}
C~Grebogi, E~Ott, and J~A Yorke.
\newblock Are three-frequency quasiperiodic orbits to be expected in typical
  nonlinear dynamical systems?
\newblock {\em Physical Review Letters}, 51(5):339, 1983.

\bibitem{nstoro6}
G.W. Luo, Y.D. Chu, Y.L. Zhang, and J.G. Zhang.
\newblock Double {N}eimark–{S}acker bifurcation and torus bifurcation of a
  class of vibratory systems with symmetrical rigid stops.
\newblock {\em Journal of Sound and Vibration}, 298(1):154--179, 2006.

\bibitem{Vibr15}
T~Bakri, Y~A Kuznetsov, and F~Verhulst.
\newblock Torus bifurcations in a mechanical system.
\newblock {\em Journal of Dynamics and Differential Equations}, 27:371--403,
  2015.

\bibitem{zhong1998torus}
Guo-Qun Zhong, Chai~Wah Wu, and Leon~O Chua.
\newblock Torus-doubling bifurcations in four mutually coupled chua's circuits.
\newblock {\em IEEE Transactions on Circuits and Systems I: Fundamental Theory
  and Applications}, 45(2):186--193, 1998.

\bibitem{kuznetsov2016simplest}
A~P Kuznetsov and Y~V Sedova.
\newblock The simplest map with three-frequency quasi-periodicity and
  quasi-periodic bifurcations.
\newblock {\em International Journal of Bifurcation and Chaos}, 26(08):1630019,
  2016.

\bibitem{Kom16}
M.~Komuro, K.~Kamiyama, T.~Endo, and K.~Aihara.
\newblock Quasi-periodic bifurcations of higher-dimensional tori.
\newblock {\em International Journal of Bifurcation and Chaos}, 26(07):1630016,
  2016.

\bibitem{Hammel94}
J.F Heagy and S.M Hammel.
\newblock The birth of strange nonchaotic attractors.
\newblock {\em Physica D: Nonlinear Phenomena}, 70(1):140--153, 1994.

\bibitem{FEU95}
U.Feudel, J.~Kurths, and A.~S. Pikovsky.
\newblock Strange non-chaotic attractor in a quasiperiodically forced circle
  map.
\newblock {\em Physica D: Nonlinear Phenomena}, 88(3):176--186, 1995.

\bibitem{SNAbook06}
U.~Feudel, S.~Kuznetsov, and A.~Pikovsky.
\newblock {\em Strange Nonchaotic Attractors}.
\newblock World Scientific, 2006.

\bibitem{kuznetsov2000critical}
S~P Kuznetsov, E~Neumann, A~Pikovsky, and I~R Sataev.
\newblock Critical point of tori collision in quasiperiodically forced systems.
\newblock {\em Physical Review E}, 62(2):1995, 2000.

\bibitem{Ga94}
J.A.C. Gallas.
\newblock Dissecting shrimps: Results for some one-dimensional physical
  systems.
\newblock {\em Physica A}, 202:196–223, 1994.

\bibitem{Gardini12}
{L. Gardini} and {I. Sushko}.
\newblock Doubling bifurcation of a closed invariant curve in 3d maps.
\newblock {\em ESAIM: Proc.}, 36:180--188, 2012.

\bibitem{Kara21}
E.~Karatetskaia, A.~Shykhmamedov, and A.~Kazakov.
\newblock {S}hilnikov attractors in three-dimensional orientation-reversing
  maps.
\newblock {\em Chaos: An Interdisciplinary Journal of Nonlinear Science},
  31(1):011102, 2021.

\bibitem{MuMcSi21}
S.S. Muni, R.I. McLachlan, and D.J.W. Simpson.
\newblock Homoclinic tangencies with infinitely many asymptotically stable
  single-round periodic solutions.
\newblock {\em Discrete Contin. Dyn. Syst. Ser A}, 41(8):3629--3650, 2021.

\bibitem{Richter02}
H.~Richter.
\newblock The generalized {H}\'{e}non maps: Examples for higher-dimensional
  chaos.
\newblock {\em International Journal of Bifurcation and Chaos},
  12(06):1371--1384, 2002.

\bibitem{de2011local}
S.~De, P~S Dutta, S~Banerjee, and A~R Roy.
\newblock Local and global bifurcations in three-dimensional, continuous,
  piecewise smooth maps.
\newblock {\em International Journal of Bifurcation and Chaos},
  21(06):1617--1636, 2011.

\bibitem{BCNF99}
S.~Banerjee and C.~Grebogi.
\newblock Border collision bifurcations in two-dimensional piecewise smooth
  maps.
\newblock {\em Phys. Rev. E}, 59:4052--4061, 1999.

\bibitem{Zhu08}
Z.~T. Zhusubaliyev and E~Mosekilde.
\newblock Formation and destruction of multilayered tori in coupled map
  systems.
\newblock {\em Chaos: An Interdisciplinary Journal of Nonlinear Science},
  18(3):037124, 2008.

\end{thebibliography}

\end{document}